\documentclass[10pt, a4paper]{article}

\usepackage{amsfonts, latexsym, amsmath, amssymb, fancyhdr, amsthm, graphicx, subfigure}
\usepackage[english]{babel}

\newtheorem{theorem}{Theorem}[section]
\newtheorem{lemma}{Lemma}[section]
\newtheorem{proposition}{Proposition}[section]

\title{\textbf{Stability of Limit Cycles in a Pluripotent Stem Cell Dynamics Model}}
\author{Mostafa Adimy$^{*}$, \quad Fabien Crauste$^{*}$, \quad Andrei Halanay$^{\dagger}$,\vspace{1ex}\\ \quad Mihaela Neam\c tu$^{\ddagger}$ \quad and
\quad Dumitru Opri\c s$^{\diamond}$}
\date{Year 2005}

\begin{document}

\maketitle

\begin{center}
{\large $^{*}$}\emph{Laboratoire de Math\'ematiques Appliqu\'ees, UMR 5142,} \\
\emph{Universit\'e de Pau et des Pays de l'Adour,}\\
\emph{Avenue de l'universit\'e, 64000 Pau, France.}\\
\emph{ANUBIS project, INRIA--Futurs}\\
\emph{Emails: mostafa.adimy@univ-pau.fr, fabien.crauste@univ-pau.fr}\\
\quad\\
{\large $^{\dagger}$}\emph{Department of Mathematics 1, University Politehnica of
Bucharest,}\\\emph{Splaiul Independen\c tei 313, RO-060042, Bucharest, Romania.}\\\emph{ Email:
halanay@vectron.mathem.pub.ro}\\
\quad\\
{\large $^{\ddagger}$}\emph{Faculty of Economics, I.N. Pestalozzi 16, West University of Timi\c
soara,}\\\emph{RO-300115, Timi\c soara, Romania.}\\\emph{ Email:
mihaela.neamtu@fse.uvt.ro}\\
\quad\\
{\large $^{\diamond}$}\emph{Department of Applied Mathematics, Faculty of Mathematics,}\\\emph{Bd.
V. Parvan 4, West University of Timi\c soara,}\\\emph{ RO-300223, Timi\c soara, Romania.}\\\emph{
Email: opris@math.uvt.ro}
\end{center}

\quad

\begin{abstract}
This paper is devoted to the study of the stability of limit cycles of a nonlinear delay
differential equation with a distributed delay. The equation arises from a model of population
dynamics describing the evolution of a pluripotent stem cells population. We study the local
asymptotic stability of the unique nontrivial equilibrium of the delay equation and we show that
its stability can be lost through a Hopf bifurcation. We then investigate the stability of the
limit cycles yielded by the bifurcation using the normal form theory and the center manifold
theorem. We illustrate our results with some numerics.
\end{abstract}

\bigskip{}


\section{Introduction}

This paper is devoted to the analysis of the nonlinear delay differential equation
\begin{equation} \label{equationx}
x^{\prime}(t)=-\big( \delta+\beta(x(t)) \big)x(t) + \frac{2}{\tau}\int_{0}^{\tau}
\beta(x(t-s))x(t-s) ds.
\end{equation}
This equation arises from a model of pluripotent hematopoietic stem cells dynamics, that is stem
cells at the root of the blood production process \cite{acr,acrsiam}. It describes the fact that
the cell density evolves according to mortality and cell division. One may stress that the cell
density considered in equation (\ref{equationx}) is in fact the density of resting cells, in
opposition to the density of proliferating cells.

The distinction between these two stages of the cell cycle is now widely accepted. We can cite, for
example, the works of Burns and Tannock \cite{burnstannock} on the existence of a resting phase ---
also called $G_0$-phase --- in the cell cycle. This phase is a quiescent stage in the cell
development, contrary to the proliferating phase which represents the active part of the cell
cycle: it is composed of the main stages of the cell development (e.g. DNA synthesis) and ends at
mitosis with the cell division. Thus, proliferating pluripotent hematopoietic stem cells are
committed to divide and give birth to two daughter cells which immediately enter the resting phase
and complete the cycle.

Mathematical models describing the dynamics of hematopoietic stem cells population have been
studied since the end of the seventies and the works of Mackey \cite{m1978,m1979}. For the reader
interested in this topic, we mention the review articles by Haurie {\it et al.} \cite{hdm98} and
Mackey {\it et al.} \cite{m2003}, and the references therein. Recently, Pujo-Menjouet {\it et al.}
\cite{pbm,pm2004} proved the existence of a Hopf bifurcation for the hematopoiesis model proposed
in \cite{m1978}, described by a nonlinear differential equation with discrete delay. However, their
results cannot be directly applied to (\ref{equationx}) because of the nature of the delay.

Delay differential equations with distributed delay have been studied by many authors. We can cite,
for example, the works in \cite{a1991,a1992,bbm2001,b1989,k1994}. However, these studies mainly
focused on stability conditions. In 2003, Liao {\it et al.} \cite{lww2003} showed the existence of
a Hopf bifurcation for a Van der Pol equation with distributed delay and studied the stability of
limit cycles, applying the normal form theory and the center manifold theorem. The characteristic
equation in \cite{lww2003} is an exponential polynomial, similar to the one obtained in
\cite{pbm,pm2004} except that the degree is higher, which makes the study easier than with equation
(\ref{equationx}). In \cite{acrsiam}, Adimy {\it et al.} obtained the existence of a Hopf
bifurcation for a nonlinear differential equation with a delay distributed according to a density,
generalizing equation (\ref{equationx}). However, these authors did not study the limit cycles of
their model.

We consider a pluripotent hematopoietic stem cells population density $x(t)$, satisfying the
nonlinear delay differential equation (\ref{equationx}).
 The constant $\delta$ accounts for natural mortality and cellular differentiation. The nonlinear
term $\beta(x(t))$ represents a rate of introduction in the proliferating phase. The last term
appears to describe the amount of cells due to cell division: cells are assumed to divide uniformly
on an interval $(0,\tau)$, with $\tau>0$, and dividing cells are in fact cells introduced in the
proliferating phase one generation earlier. The assumption on the cell division comes from the fact
that, even though only a little is known about phenomena involved in hematopoiesis, there are
strong evidences (see Bradford \emph{et al.} \cite{bradford}) indicating that cells do not divide
at the same age. The factor $2$ describes the division of each mother cell in two daughter cells.


The rate of reintroduction in the proliferating compartment $\beta$ is taken to be a monotone and
decreasing Hill function, given by
\begin{equation}\label{beta}
\beta(x)=\beta_0\frac{\theta^n}{\theta^n+x^n} \qquad \textrm{ for } x\geq0.
\end{equation}
The coefficient $\beta_0>0$ is the maximum rate of reintroduction, $\theta\geq0$ is the resting
phase population density for which the rate of re-entry $\beta$ attains its maximum rate of change
with respect to the resting phase population, and $n\geq0$ describes the sensitivity of $\beta$
with changes in the population. This function was firstly used in hematopoiesis models by Mackey
\cite{m1978} in 1978.

In \cite{m1978,m1979} Mackey gave values of the above parameters for a normal human body
production. These values are
\begin{equation}\label{parametersvalues}
\delta=0.05\textrm{ d}^{-1} \quad \textrm{ and } \quad \beta_0=1.77\textrm{ d}^{-1}.
\end{equation}
Usually, $n$ is close to $1$, but Mackey \cite{m1978,m1979} reports values of $n$ around $3$ in
abnormal situations.

\begin{figure}[!hpt]
\begin{center}
\includegraphics[width=10cm,height=7.5cm]{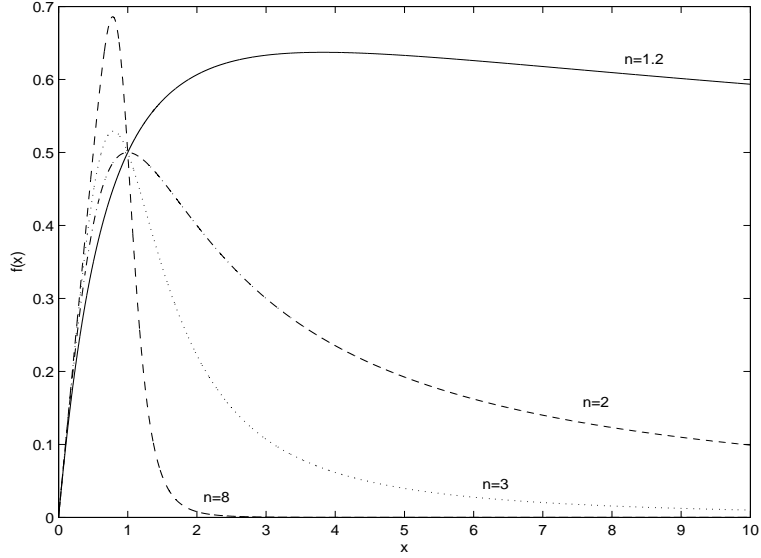}
\end{center}
\caption{Graph of the function $f$ given by (\ref{functionf}) for different values of $n$.}
\label{figf}
\end{figure}

The value of $\theta$ is usually $\theta=1.62\times 10^8\textrm{ cells/kg}.$ However, since we
shall study the qualitative behavior of the pluripotent stem cells population, the value of
$\theta$ is not really important and setting the scale change
\begin{displaymath}
x(t) \to \frac{x(t)}{\theta}
\end{displaymath}
in (\ref{equationx}), with the function $\beta$ given by (\ref{beta}), we obtain
\begin{equation}\label{eq}
\frac{dx}{dt}(t)=-\delta x(t) - \beta_0 f\left(x(t)\right) + \frac{2\beta_0}{\tau}\int^0_{-\tau}
f\left(x(t+s)\right)ds
\end{equation}
with
\begin{equation}\label{functionf}
f(x)=\frac{x}{1+x^n}, \qquad x\geq0.
\end{equation}
However, we mention that this special form of $f$ will not be used in the following, except in
computations in Section \ref{snumerical}. We only assume that $f$ is differentiable with $f(0)=0$
and, for $x>0$, $f(x)/x$ is decreasing and satisfies
\begin{displaymath}
\lim_{x\to +\infty}\frac{f(x)}{x}=0.
\end{displaymath}
An illustration of such a function is presented in Fig. \ref{figf}.

Notice that equation (\ref{eq}) has at most two equilibria: the trivial equilibrium $x\equiv 0$ and
a nontrivial positive equilibrium $x\equiv x^*$. The trivial equilibrium always exists and
corresponds to the extinction of the population.

The nontrivial equilibrium exists if and only if
\begin{equation}\label{existencexstar}
0<\delta<\beta_0f^{\prime}(0)
\end{equation}
and is then uniquely defined by
\begin{equation}\label{xstar}
\delta x^*=\beta_0f(x^*).
\end{equation}
This can be easily shown by using the fact that the function $f(x)/x$ is decreasing on
$[0,+\infty)$.

Our aim in this work is to show that the unique nontrivial equilibrium of equation
(\ref{equationx}) undergoes, in a particular case, a unique Hopf bifurcation and to show the
stability of the limit cycles following the approach in \cite{hkw1981,s1989}.

The paper is organized as follows. In Section \ref{shopf}, we establish some local stability
results for the unique nontrivial equilibrium of (\ref{equationx}) and prove that it undergoes a
Hopf bifurcation. We study the stability of the limit cycles obtained at the bifurcation in Section
\ref{sstability}. Our results are illustrated numerically in Section \ref{snumerical}. We conclude
with a discussion in Section \ref{sdiscussion}.


\section{Local Stability and Hopf Bifurcation Analysis}\label{shopf}

Part of the study presented in this section has been previously performed by Adimy {\it et al.} in
\cite{acr}. However, for the reader convenience and to preserve the coherence of the present work,
we detail the asymptotic behavior study of the nontrivial equilibrium $x\equiv x^*$ of equation
(\ref{eq}), defined by (\ref{xstar}). We first concentrate on the local asymptotic stability of
this equilibrium. Then, we will show that it undergoes a Hopf bifurcation for some critical value
of the time delay.

We assume that (\ref{existencexstar}) holds in order to ensure the existence of $x^*$, that is
\begin{displaymath}
0<\delta<\beta_0f^{\prime}(0).
\end{displaymath}
The linearization of equation (\ref{eq}) around $x^*$ leads to the characteristic equation
\begin{equation}\label{ce}
\Delta(\lambda,\tau):=\lambda+\delta+\beta_0\beta_1-\frac{2\beta_0\beta_1}
{\tau}\int^{0}_{-\tau}e^{\lambda s}ds=0,
\end{equation}
where we have set
$$
\beta_1:=f^{\prime}(x^*).
$$
We recall that the nontrivial equilibrium $x^*$ is locally asymptotically stable if and only if all
eigenvalues of (\ref{ce}) have negative real parts.

The function $f$ in (\ref{eq}) is not necessarily monotone so $\beta_1$ may be either positive or
negative. We first study the case $\beta_1\geq0$.

Consider $\Delta(\lambda,\tau)$ as a function of real $\lambda$. Then $\Delta$ is differentiable
with respect to $\lambda$ and
\begin{equation}\label{ced}
\Delta_{\lambda}(\lambda,\tau):=\frac{\partial\Delta}{\partial\lambda}(\lambda,\tau)=1-\frac{2\beta_0\beta_1}{\tau}
\int^{0}_{-\tau}se^{\lambda s}ds.
\end{equation}
Since $\displaystyle\int^{0}_{-\tau}se^{\lambda s}ds<0,$ we deduce that
$\Delta_{\lambda}(\lambda,\tau)>0$. Moreover, one can easily check that
\begin{displaymath}
\lim_{\lambda\to-\infty}\Delta(\lambda,\tau)=-\infty \qquad \textrm{ and } \qquad
\lim_{\lambda\to+\infty}\Delta(\lambda,\tau)=+\infty.
\end{displaymath}
Consequently, $\Delta(\lambda,\tau)$ has a unique real eigenvalue, namely $\lambda_0$. From
(\ref{ce}), we obtain
\begin{displaymath}
\Delta(0,\tau)=\delta-\beta_0\beta_1.
\end{displaymath}
Writing
\begin{displaymath}
\beta_1=x^*\left(\frac{f(x)}{x}\right)^{\prime}\bigg|_{x=x^*}+\frac{f(x^*)}{x^*}
=x^*\left(\frac{f(x)}{x}\right)^{\prime}\bigg|_{x=x^*}+\frac{\delta}{\beta_0},
\end{displaymath}
we deduce
\begin{displaymath}
\Delta(0,\tau)=-\beta_0x^*\left(\frac{f(x)}{x}\right)^{\prime}\bigg|_{x=x^*}>0.
\end{displaymath}
Hence, $\lambda_0$ is strictly negative.

Let $\lambda=\nu+i\omega$ be an eigenvalue of (\ref{ce}) and assume that $\nu>\lambda_0$.
Considering the real part of (\ref{ce}), we obtain
\begin{displaymath}
\nu+\delta+\beta_0\beta_1-\frac{2\beta_0\beta_1} {\tau}\int^{0}_{-\tau}e^{\nu s}\cos(\omega s)ds=0.
\end{displaymath}
Therefore,
\begin{displaymath}
\nu-\lambda_0=\frac{2\beta_0\beta_1} {\tau}\int^{0}_{-\tau}\left[e^{\nu s}\cos(\omega
s)-e^{\lambda_0 s}\right]ds<0,
\end{displaymath}
which gives a contradiction. We deduce that every eigenvalue of (\ref{ce}) has negative real part.
It follows that the nontrivial positive equilibrium $x^*$ is locally asymptotically stable. This
result is summed up in the following proposition.

\begin{proposition}
Assume that $\beta_1\geq0$. Then the nontrivial equilibrium $x\equiv x^*$ of (\ref{eq}) is locally
asymptotically stable for all $\tau\geq0$.
\end{proposition}

If the function $f$ is given by (\ref{functionf}) then (\ref{existencexstar}) is equivalent to
$0<\delta<\beta_0$ and $x^*=(\beta_0/\delta-1)^{1/n}$. The condition $\beta_1\geq0$ reduces to
\begin{displaymath}
n\frac{\beta_0-\delta}{\beta_0}\leq1.
\end{displaymath}

We assume now that $\beta_1<0$. We are going to show that the equilibrium $x^*$ undergoes a Hopf
bifurcation. To that aim, we look for the existence of purely imaginary roots of (\ref{ce}).

We first check that $x^*$ is locally asymptotically stable when $\tau=0$. In this case, the
characteristic equation (\ref{ce}) reduces to
\begin{displaymath}
\lambda+\delta-\beta_0\beta_1=0,
\end{displaymath}
so
\begin{displaymath}
\lambda=-\delta+\beta_0\beta_1<0.
\end{displaymath}
We have the following lemma.

\begin{lemma}\label{lemma}
Assume that $\beta_1<0$. Then the nontrivial equilibrium $x\equiv x^*$ of (\ref{eq}) is locally
asymptotically stable when $\tau=0$.
\end{lemma}

Let $\lambda=i\omega$, $\omega\in\mathbb{R}$, be a purely imaginary eigenvalue of (\ref{ce}). One
can check that
\begin{displaymath}
\Delta(-i\omega,\tau)=0,
\end{displaymath}
so we only look for positive $\omega$. Moreover, $\omega\neq0$ since
\begin{displaymath}
\Delta(0,\tau)=\delta-\beta_0\beta_1>0.
\end{displaymath}
Thus, let us assume that $\omega>0$ and $\tau>0$ satisfy $\Delta(i\omega,\tau)=0$. Separating real
and imaginary parts of $\Delta(i\omega,\tau)$, we obtain
\begin{equation}\label{system1}
\left\{\begin{array}{rcl}
\delta+\beta_0\beta_1-\displaystyle\frac{2\beta_0\beta_1}{\omega\tau}\sin(\omega\tau)&=&0,\vspace{1ex}\\
\omega+\displaystyle\frac{2\beta_0\beta_1}{\omega\tau}(1-\cos(\omega\tau))&=&0.
\end{array}\right.
\end{equation}
We set
\begin{displaymath}
h(x)=\frac{\sin(x)}{x}, \qquad x>0.
\end{displaymath}
Then system (\ref{system1}) can be written
\begin{equation}\label{system2}
\left\{\begin{array}{rcl}
h(\omega\tau)&=&\displaystyle\frac{\delta+\beta_0\beta_1}{2\beta_0\beta_1},\vspace{1ex}\\
\displaystyle\frac{\cos(\omega\tau)-1}{(\omega\tau)^2}&=&\displaystyle\frac{1}{2\beta_0\beta_1\tau}.
\end{array}\right.
\end{equation}
Since $\beta_1<0$ and $\delta>0$, then
\begin{displaymath}
\frac{\delta+\beta_0\beta_1}{2\beta_0\beta_1}<\frac{1}{2}.
\end{displaymath}
Let
\begin{displaymath}
x_0=\min\left\{x>0;\ x=\tan(x) \textrm{ and } h(x)>0 \right\} \simeq 7.725,
\end{displaymath}
and assume that
\begin{equation}\label{hyph}
h(x_0)<\frac{\delta+\beta_0\beta_1}{2\beta_0\beta_1}.
\end{equation}
One can check that $h(x_0)\simeq 0.1284$. Then (\ref{hyph}) is equivalent to
\begin{displaymath}
\beta_1<-\frac{\delta}{\beta_0(1-2h(x_0))}.
\end{displaymath}
On the interval $[0,\pi]$, the function $h$ is strictly decreasing and nonnegative, with $0\leq
h(x)\leq1$. Moreover, for $x\geq\pi$, $h(x)\leq h(x_0)$. Consequently, the equation
\begin{displaymath}
h(x)=\frac{\delta+\beta_0\beta_1}{2\beta_0\beta_1}
\end{displaymath}
has a unique solution, denoted $x_c$, which belongs to the interval $(0,\pi)$. We then set
\begin{displaymath}
\tau_c=\frac{x_c^2}{2\beta_0\beta_1(\cos(x_c)-1)}
\end{displaymath}
and
\begin{displaymath}
\omega_c=\frac{x_c}{\tau_c}.
\end{displaymath}
Therefore, $(\omega_c,\tau_c)$ is the unique solution of (\ref{system2}) and $\pm i\omega_c$ are
purely imaginary eigenvalues of (\ref{ce}) for $\tau=\tau_c$.

In order to show that $x^*$ undergoes a Hopf bifurcation for $\tau=\tau_c$, we have to prove that
$\pm i\omega_c$ are simple eigenvalues of $\Delta(\cdot,\tau_c)$ and satisfy the transversality
condition
\begin{displaymath}
\frac{d \textrm{Re}(\lambda)}{d\tau}\bigg|_{\tau=\tau_c}>0.
\end{displaymath}

We first check that $\pm i\omega_c$ are simple eigenvalues of (\ref{ce}). Using (\ref{ced}), one
can see that $i\omega_c$ is simple if
\begin{equation}\label{eq0}
\textrm{Re}\left(\Delta_{\lambda}(i\omega_c,\tau_c)\right)=1+\frac{2\beta_0\beta_1}{\tau_c}\
\frac{\omega_c\tau_c\sin(\omega_c\tau_c)+\cos(\omega_c\tau_c)-1}{\omega_c^2}\neq0
\end{equation}
or
\begin{equation}\label{eq1}
\textrm{Im}\left(\Delta_{\lambda}(i\omega_c,\tau_c)\right)=\frac{2\beta_0\beta_1}{\tau_c}\
\frac{\omega_c\tau_c\cos(\omega_c\tau_c)-\sin(\omega_c\tau_c)}{\omega_c^2}\neq0.
\end{equation}
We are going to show that, in fact, these two conditions are satisfied.

\begin{lemma}\label{lemma2}
Assume that $\beta_1<0$ and (\ref{hyph}) holds. Let $(\omega_c,\tau_c)$ be the unique solution of
(\ref{system2}), with $\omega_c\tau_c\in(0,\pi)$. Then
\begin{equation}\label{eq1bis}
\textrm{Re}\left(\Delta_{\lambda}(i\omega_c,\tau_c)\right)>0 \qquad \textrm{ and } \qquad
\textrm{Im}\left(\Delta_{\lambda}(i\omega_c,\tau_c)\right)>0.
\end{equation}
In particularly, $\pm i\omega_c$ are simple eigenvalues of (\ref{ce}) for $\tau=\tau_c$.
\end{lemma}

\begin{proof}
First, one can check, using (\ref{system2}) and (\ref{eq0}), that
\begin{displaymath}
\textrm{Re}\left(\displaystyle\Delta_{\lambda}(i\omega_c,\tau_c)\right)=2+(\delta+\beta_0\beta_1)\tau_c.
\end{displaymath}
Since
\begin{displaymath}
\frac{\cos(x)-1}{x^2}=-\frac{h^2(x)}{1+\cos(x)}, \qquad \textrm{ for } x>0,
\end{displaymath}
then, from (\ref{system2}) it follows that
\begin{displaymath}
1+\cos(x)=-(\delta+\beta_0\beta_1)\tau_c.
\end{displaymath}
Consequently,
\begin{displaymath}
\textrm{Re}\left(\displaystyle\Delta_{\lambda}(i\omega_c,\tau_c)\right)=1-\cos(\omega_c\tau_c)>0.
\end{displaymath}

Secondly, since $x\cos(x)<\sin(x)$ for $x\in(0,\pi)$ and $x_c=\omega_c\tau_c\in(0,\pi)$, then from
(\ref{eq1}) we obtain
\begin{displaymath}
\textrm{Im}\left(\Delta_{\lambda}(i\omega_c,\tau_c)\right)>0.
\end{displaymath}
This concludes the proof.
\end{proof}

Consider now a branch of eigenvalues $\lambda(\tau)=\nu(\tau)+i\omega(\tau)$ of (\ref{ce}) such
that $\nu(\tau_c)=0$ and $\omega(\tau_c)=\omega_c$. Separating real and imaginary parts in
(\ref{ce}) we obtain
\begin{displaymath}
\left\{\begin{array}{rcl} \nu(\tau)+\delta+\beta_0\beta_1-\displaystyle\frac{2\beta_0\beta_1}{\tau}
\int^0_{-\tau}e^{\nu(\tau) s}\cos(\omega(\tau) s)ds&=&0,\vspace{1ex}\\
\omega(\tau)-\displaystyle\frac{2\beta_0\beta_1}{\tau}\int^0_{-\tau}e^{\nu(\tau)
s}\sin(\omega(\tau) s)ds&=&0.
\end{array}\right.
\end{displaymath}
Then, by differentiating each of the above equalities with respect to $\tau$, we get, for
$\tau=\tau_c$,
\begin{equation}\label{muprimezero}
\begin{array}{l}
\textrm{Re}\left(\displaystyle\Delta_{\lambda}(i\omega_c,\tau_c)\right)\nu^{\prime}(\tau_c)\vspace{1ex}\\
=\textrm{Im}\left(\displaystyle\Delta_{\lambda}(i\omega_c,\tau_c)\right)\omega^{\prime}(\tau_c)
+\displaystyle\frac{2\beta_0\beta_1}{\tau_c}\left(\cos(x_c)-\frac{\sin(x_c)}{x_c}\right)
\end{array}
\end{equation}
and
\begin{equation}\label{omegaprimezero}
\begin{array}{l}
\textrm{Re}\left(\displaystyle\Delta_{\lambda}(i\omega_c,\tau_c)\right)\omega^{\prime}(\tau_c)\vspace{1ex}\\
=-\textrm{Im}\left(\displaystyle\Delta_{\lambda}(i\omega_c,\tau_c)\right)\nu^{\prime}(\tau_c)
+\displaystyle\frac{2\beta_0\beta_1}{\tau_c}\left(\frac{1-\cos(x_c)}{x_c}-\sin(x_c)\right).
\end{array}
\end{equation}
Using (\ref{eq1bis}), (\ref{muprimezero}) and (\ref{omegaprimezero}), we can see that
$\nu^{\prime}(\tau_c)$ satisfies
\begin{displaymath}
\begin{array}{l}
\left[\textrm{Im}\left(\displaystyle\Delta_{\lambda}(i\omega_c,\tau_c)\right)^2
+\textrm{Re}\left(\displaystyle\Delta_{\lambda}(i\omega_c,\tau_c)\right)^2\right]
\nu^{\prime}(\tau_c)\vspace{1ex}\\
=\displaystyle\frac{2\beta_0\beta_1}{\tau_c}\left[ \left(\frac{1-\cos(x_c)}{x_c}-\sin(x_c)\right)
\textrm{Im}\left(\displaystyle\Delta_{\lambda}(i\omega_c,\tau_c)\right)\right.\vspace{1ex}\\
\left.\qquad\qquad\qquad+\left(\cos(x_c)-\displaystyle\frac{\sin(x_c)}{x_c}\right)
\textrm{Re}\left(\displaystyle\Delta_{\lambda}(i\omega_c,\tau_c)\right) \right].
\end{array}
\end{displaymath}
Using the definitions in (\ref{eq0}) and (\ref{eq1}), simple computations give
\begin{displaymath}
\begin{array}{rcl}
\left(\displaystyle\frac{1-\cos(x_c)}{x_c}-\sin(x_c)\right)
\textrm{Im}\left(\displaystyle\Delta_{\lambda}(i\omega_c,\tau_c)\right)&&\vspace{1ex}\\
+\left(\cos(x_c)-\displaystyle\frac{\sin(x_c)}{x_c}\right)
\textrm{Re}\left(\displaystyle\Delta_{\lambda}(i\omega_c,\tau_c)\right)&=&\cos(x_c)-\displaystyle\frac{\sin(x_c)}{x_c}.
\end{array}
\end{displaymath}
Hence,
\begin{displaymath}
\left|\displaystyle\Delta_{\lambda}(i\omega_c,\tau_c)\right|^2\nu^{\prime}(\tau_c)=
\displaystyle\frac{2\beta_0\beta_1}{\tau_c}\ \displaystyle\frac{x_c\cos(x_c)-\sin(x_c)}{x_c}>0.
\end{displaymath}
It follows that
\begin{equation}\label{muprime}
\nu^{\prime}(\tau_c)>0.
\end{equation}

To conclude, when $\tau=\tau_c$, the characteristic equation $\Delta(\lambda,\tau)$ has a unique
pair of purely imaginary simple eigenvalues satisfying
$(d\textrm{Re}(\lambda)/d\tau)(\tau=\tau_c)>0$. Consequently, a Hopf bifurcation occurs at $x^*$
when $\tau=\tau_c$. Moreover, applying Rouch\'e's Theorem with Lemma \ref{lemma}, we easily check
that every eigenvalue of $\Delta(\lambda,\tau)$, with $\tau<\tau_c$, has negative real part. It
follows that $x^*$ is locally asymptotically stable for $0\leq\tau<\tau_c$. These results are
summed up in the following theorem.

\begin{theorem}\label{theoremhopfbifurcation}
Assume that $\beta_1<0$ and (\ref{hyph}) holds. Then there exists a unique value $\tau_c>0$ of the
time delay such that the equilibrium $x\equiv x^*$ is locally asymptotically stable when
$\tau\in[0,\tau_c)$ and becomes unstable when $\tau=\tau_c$ throughout a Hopf bifurcation. In
particularly, periodic solutions appear for equation (\ref{eq}) when $\tau=\tau_c$.
\end{theorem}

As an example, one can check that when $f$ is given by (\ref{functionf}) the assumptions in Theorem
\ref{theoremhopfbifurcation} are equivalent to
\begin{displaymath}
n>\frac{2(1-h(x_0))}{1-2h(x_0)}\ \frac{\beta_0}{\beta_0-\delta} \simeq 2.35
\frac{\beta_0}{\beta_0-\delta}.
\end{displaymath}
In particularly, these conditions are satisfied when $\beta_0$ and $\delta$ are given by
(\ref{parametersvalues}) and $n\geq 2.42$.


The existence of a Hopf bifurcation in Theorem \ref{theoremhopfbifurcation} leads to the existence
of a limit cycle when the bifurcation occurs. In the next section, we focus on the stability of
this limit cycle.

\section{Stability of Limit Cycles}\label{sstability}

We study now the stability of the limit cycle yielded by Theorem \ref{theoremhopfbifurcation}. We
follow the approach used in \cite{hkw1981,s1989}. This involves the description of a center
manifold and subsequently the study of the normal form given by the restriction of the flow to this
center manifold. The stability of the limit cycle will be decided by the sign of the first Lyapunov
coefficient $l_1(0)$.

For general properties concerning delay equations and the theory of central manifolds for these
equations, see \cite{h1966}. For the existence and various properties of center manifolds we refer
to \cite{c1981,c1971,cm1977,fm1995,h1977,mm1976}. Also, in \cite{fm1995} and \cite{sr2003}, a
rigorous treatment of the operators involved in this approach is to be found. A similar problem is
considered in \cite{lww2003}.

Define, for $t\geq0$,
\begin{displaymath}
y(t)= x(t)-x^*,  \qquad \mu = \tau - \tau_c
\end{displaymath}
with $x^*$ the nontrivial equilibrium of (\ref{eq}) that bifurcates into a limit cycle for the
critical value $\tau=\tau_c$ (see Theorem \ref{theoremhopfbifurcation}). The equilibrium $x^*$ is
defined by (\ref{existencexstar}) and (\ref{xstar}). Equation (\ref{eq}) turns into
\begin{equation}\label{eqy}
y^{\prime}(t)=-\delta(y(t)+x^*)-\beta_0f(y(t)+x^*)+ \frac{2\beta_0}{\mu+\tau_c}\int^0_{-\mu
-\tau_c} f\left(y(t+s)+x^*\right)ds.
\end{equation}
Thanks to this formulation, we now concentrate on the trivial equilibrium $y\equiv 0$ of
(\ref{eqy}) which bifurcates when $\mu=0$.

For an interval $I\subset\mathbb{R}$, denote $C(I,\mathbb{K}) = \{f:I\to\mathbb{K}, f \textrm{
continuous } \}$ where $\mathbb{K}=\mathbb{R}$ or $\mathbb{C}$. When $I=[-\mu-\tau_c,0]$, we set
\begin{displaymath}
C_{\mu} := C([-\mu-\tau_c,0],\mathbb{K}).
\end{displaymath}

Considering, for $t\geq0$, the function $y_t:[-\mu-\tau_c,0]\to\mathbb{K}$ defined by
$y_t(s)=y(t+s)$, we can reformulate equation (\ref{eqy}) as the following abstract functionnal
differential equation
\begin{equation}\label{aeqy}
\frac{d}{dt}y(t)=G_{\mu}(y_t), \qquad t\geq0,
\end{equation}
where, for $\varphi\in C_{\mu}$,
\begin{displaymath}
G_{\mu}(\varphi)=-\delta\left[\varphi(0)+x^*\right]-\beta_0f\left(\varphi(0)+x^*\right)
+\displaystyle\frac{2\beta_0}{\mu+\tau_c}\int^0_{-\mu -\tau_c} f\left(\varphi(s)+x^*\right)ds.
\end{displaymath}
Assume that $f$ is $\mathcal{C}^4$ on $[0,+\infty)$ (remark that $G_{\mu}$ is then
$\mathcal{C}^4(C_{\mu},\mathbb{R})$).

Consider the linearized equation of (\ref{aeqy}), corresponding to the Fr\'echet derivative
$D_{\varphi}G_{\mu}(0) := L_{\mu}$, given by
\begin{equation}\label{leqy}
\frac{d}{dt}z(t)=L_{\mu}z_t, \qquad t\geq0.
\end{equation}
In fact, $L_{\mu}$ is given explicitly by
\begin{equation}\label{L}
L_{\mu}\varphi=-c_1\varphi(0)+c_2(\mu)\displaystyle\int^0_{-\mu -\tau_c} \varphi(\theta)d\theta,
\qquad \varphi\in C_{\mu},
\end{equation}
where
\begin{equation}\label{parameters1}
c_1 := \delta+ \beta_0\beta_1, \qquad c_2(\mu) := \displaystyle\frac{2\beta_0\beta_1}{\tau_c +
\mu}, \qquad \beta_1:=f'(x^*).
\end{equation}

Setting
\begin{displaymath}
F_{\mu} := G_{\mu}-L_{\mu},
\end{displaymath}
equation (\ref{aeqy}) becomes
\begin{equation}\label{aeqy2}
\frac{d}{dt}y(t)=L_{\mu}y_t+F_{\mu}(y_t), \qquad t\geq0,
\end{equation}
with $F_{\mu}(0)=0$ and $D_{\varphi}F_{\mu}(0)=0$.

In order to develop a normal form associated to equation (\ref{leqy}), we write this latter as an
abstract ordinary differential equation.

First, we know from \cite{h1977} that the linear equation (\ref{leqy}) gives a $C_0$-semigroup
$(T(t))_{t\geq0}$ on $C_{\mu}$, with generator $A_{\mu}$ defined by
\begin{displaymath}
\left\{\begin{array}{rcl} \mathcal{D}(A_{\mu})&=&\left\{\varphi\in C^1([-\mu-\tau_c,0],\mathbb{R});
\
\varphi^{\prime}(0)=L_{\mu}\varphi\right\},\vspace{1ex}\\
A_{\mu}\varphi&=&\varphi^{\prime}, \qquad \varphi\in\mathcal{D}(A_{\mu}).
\end{array}\right.
\end{displaymath}
To write equation (\ref{aeqy2}) as an ODE we need to extend the problem (\ref{leqy}) to the Banach
space $\widetilde{C}_{\mu} := C_{\mu}\oplus\langle X_0\rangle$, where
\begin{displaymath}
\langle X_0\rangle = \left\{ X_0c; \ c\in\mathbb{R} \textrm{ and } (X_0c)(\theta)=X_0(\theta)c
\right\}
\end{displaymath}
and $X_0$ denotes the function defined on $[-\mu-\tau_c,0]$ by
\begin{displaymath}
X_0(\theta)=\left\{\begin{array}{ll} 0,& \quad\textrm{ if }
-\mu-\tau_c\leq\theta<0,\vspace{1ex}\\
1,& \quad\textrm{ if } \theta=0.
\end{array}\right.
\end{displaymath}
Adimy proved in \cite{a1993} that this extension determines a Hille-Yosida operator. This result is
recalled in the next lemma.

\begin{lemma}
The continuous extension $\widetilde{A}_{\mu}$ of the operator $A_{\mu}$ defined on
$\widetilde{C}_{\mu}$ by
\begin{equation}\label{At}
\begin{array}{rcl}
\mathcal{D}(\widetilde{A}_{\mu})&=& C^1([-\mu-\tau_c,0],\mathbb{R}),\vspace{1ex}\\
\widetilde{A}_{\mu}\varphi&=&\varphi^{\prime}+X_0\left(L_{\mu}\varphi-\varphi^{\prime}(0)\right),
\qquad \varphi\in\mathcal{D}(\widetilde{A}_{\mu}),
\end{array}
\end{equation}
is a Hille-Yosida operator on $\widetilde{C}_{\mu}$; that is: there exists $\omega_0\in\mathbb{R}$
such that $(\omega_0,+\infty)\subset\rho(\widetilde{A}_{\mu})$ and
\begin{displaymath}
\sup\left\{(\lambda-\omega_0)^n\|(\lambda I-\widetilde{A}_{\mu})^{-n}\|, n\in\mathbb{N},
\lambda>\omega_0\right\}<\infty.
\end{displaymath}
\end{lemma}

It follows that if $y$ is a solution of (\ref{aeqy2}) on $[0,T]$, $T>0$, with an initial condition
$\varphi\in C_{\mu}$ on the interval $[-\mu-\tau_c,0]$, then the function $t\in[0,T]\mapsto y_t\in
C_{\mu}$ satisfies
\begin{equation}\label{aeqy3}
\left\{\begin{array}{rcl} \displaystyle\frac{d}{dt}y_t&=&\widetilde{A}_{\mu}y_t+X_0F_{\mu}(y_t),
\qquad
t\in[0,T],\vspace{1ex}\\
y_0&=&\varphi.
\end{array}\right.
\end{equation}
Conversely, if there exists a function $t\in[0,T]\mapsto u(t)\in C_{\mu}$ such that
\begin{equation}\label{equ}
\left\{\begin{array}{rcl} \displaystyle\frac{du}{dt}(t)&=&\widetilde{A}_{\mu}u(t)+X_0F_{\mu}(u(t)),
\qquad
t\in[0,T],\vspace{1ex}\\
u(0)&=&\varphi,
\end{array}\right.
\end{equation}
then $u(t)=y_t$, $t\in[0,T]$, where
\begin{displaymath}
y(t)=\left\{\begin{array}{ll}
u(t)(0),&\qquad\textrm{ if } t\in[0,T],\vspace{1ex}\\
\varphi(t),&\qquad\textrm{ if } t\in[-\mu-\tau_c,0],
\end{array}\right.
\end{displaymath}
and $y$ is a solution of (\ref{aeqy2}). This handles in particularly the problems arising from the
fact that $\widetilde{A}_{\mu}$ does not preserve the space of $\mathcal{C}^1$-functions.

Thanks to results by Arendt \cite{a1987} and Da Prato and Sinestrari \cite{ds1987}, the ODE
(\ref{equ}) is well-posed for initial conditions in
$\overline{\mathcal{D}(\widetilde{A}_{\mu})}=C_{\mu}$.


Now we can reformulate the problem (\ref{aeqy2}) as the abstract ODE (\ref{aeqy3}).

Another important step towards the description of a center manifold is the definition of a bilinear
form related to the equation (\ref{leqy}).

From now on, we set $\mathbb{K}=\mathbb{C}$. For $\varphi\in C_{\mu}$ and $\psi \in
C_{\mu}^*:=C([0,\mu+\tau_c],\mathbb{C})$, define according to \cite{fm1995} or \cite{h1977},
\begin{equation}\label{ps}
\langle\psi,\varphi\rangle = \overline{\psi(0)}\varphi(0) - \displaystyle\int^0_{-\mu-\tau_c}
\left(\int^s_0 \overline{\psi(\xi-s)}\varphi(\xi)d\xi \right) d\eta (s)
\end{equation}
where $d\eta (s) = c_2(\mu)ds-c_1X_0(s)$. Thus (\ref{ps}) becomes
\begin{equation}\label{ps2}
\langle\psi,\varphi\rangle = \overline{\psi (0)} \varphi(0) - c_2 (\mu)
\displaystyle\int^0_{-\mu-\tau_c} \left(\int^s_0 \overline{\psi(\xi-s)}\varphi(\xi)d\xi \right) ds.
\end{equation}
We build a natural extension of this bilinear form to the space
$\widetilde{C}_{\mu}^*\times\widetilde{C}_{\mu}$ where
\begin{displaymath}
\widetilde{C}_{\mu}^* = C_{\mu}^*\oplus\langle X_0^*\rangle
\end{displaymath}
and
\begin{displaymath}
\langle X_0^*\rangle = \left\{ X_0^*c; \ c\in\mathbb{C} \textrm{ and }
(X_0^*c)(\theta)=X_0^*(\theta)c \right\}
\end{displaymath}
with $X_0^*$ the function defined on $[0,\mu+\tau_c]$ by
\begin{displaymath}
X_0^*(\theta)=\left\{\begin{array}{ll} 0,& \quad\textrm{ if }
0<\theta\leq \mu+\tau_c,\vspace{1ex}\\
1,& \quad\textrm{ if } \theta=0.
\end{array}\right.
\end{displaymath}
We obtain, for $\psi\in C_{\mu}^*$, $\varphi\in C_{\mu}$ and $a,c\in\mathbb{C}$,
\begin{displaymath}
\langle \psi+X_0^*a,\varphi+X_0c \rangle=\langle \psi,\varphi \rangle + \overline{a}c.
\end{displaymath}

With respect to this bilinear form, we define the adjoint of the operator $\widetilde{A}_{\mu}$,
denoted $\widetilde{A}_{\mu}^*$, and its domain $\mathcal{D}(\widetilde{A}_{\mu}^*)$. It satisfies,
for $\varphi\in\mathcal{D}(A_{\mu})=C^1([-\mu-\tau_c,0],\mathbb{C})$ and
$\psi\in\mathcal{D}(\widetilde{A}_{\mu}^*)$,
\begin{displaymath}
\langle \psi,\widetilde{A}_{\mu}\varphi \rangle = \langle \widetilde{A}_{\mu}^*\psi,\varphi
\rangle.
\end{displaymath}
From (\ref{At}) and (\ref{ps2}), we obtain
\begin{displaymath}
\langle \psi,\widetilde{A}_{\mu}\varphi \rangle = \overline{\psi (0)}L_{\mu}\varphi-c_2 (\mu)
\displaystyle\int^0_{-\mu-\tau_c} \left(\int^s_0 \overline{\psi(\xi-s)}\varphi^{\prime}(\xi)d\xi
\right) ds.
\end{displaymath}
Using an integration by parts and (\ref{L}), we deduce
\begin{displaymath}
\setlength\arraycolsep{2pt}
\begin{array}{rcl}
\langle \psi,\widetilde{A}_{\mu}\varphi \rangle &=& \overline{\psi(0)}
\left[-c_1\varphi(0)+c_2(\mu)\displaystyle\int^0_{-\mu-\tau_c}\varphi(\theta)d\theta\right]\vspace{1ex}\\
&&-c_2(\mu)\displaystyle\int^0_{-\mu-\tau_c} \left(
\overline{\psi(0)}\varphi(s)-\overline{\psi(-s)}\varphi(0)-
\int^s_0\overline{\psi^{\prime}(\xi-s)}\varphi(\xi)d\xi
\right)ds,\vspace{1ex}\\
&=&\left[-c_1\overline{\psi(0)}+c_2(\mu)\displaystyle\int_0^{\mu+\tau_c}\overline{\psi(\theta)}d\theta\right]\varphi(0)\vspace{1ex}\\
&&+c_2(\mu)\displaystyle\int^0_{-\mu-\tau_c} \left(
\int^s_0\overline{\psi^{\prime}(\xi-s)}\varphi(\xi)d\xi
\right)ds,\vspace{1ex}\\
&=&\langle \widetilde{A}_{\mu}^*\psi,\varphi \rangle,
\end{array}
\end{displaymath}
where
\begin{displaymath}
\left\{\begin{array}{rcl}
\mathcal{D}(\widetilde{A}_{\mu}^*)&=&C^1([0,\mu+\tau_c],\mathbb{R}),\vspace{1ex}\\
\widetilde{A}_{\mu}^*\psi&=&-\psi^{\prime}+X_0^*\left[c_2 (\mu) \displaystyle\int^{\mu + \tau_c}_0
\psi(s)ds-c_1\psi(0)+\psi^{\prime}(0)\right].
\end{array}\right.
\end{displaymath}

We consider now the purely imaginary eigenvalues of (\ref{ce}) denoted $\pm i\omega_c$, with
$\omega_c>0$, which exist when $\tau=\tau_c$, that means when the bifurcation occurs (see Section
\ref{shopf} and, in particularly, Theorem \ref{theoremhopfbifurcation}). From the definition of the
characteristic equation in (\ref{ce}) and the notations introduced in (\ref{parameters1}), we have
\begin{displaymath}
\Delta(i\omega_c, \tau_c)=i\omega_c+c_1-c_2(0)\int^{0}_{-\tau_c}e^{i\omega_c s}ds=0.
\end{displaymath}
It follows that
\begin{equation}\label{prop}
c_2 (0) \frac{(e^{-i\omega_c\tau_c}-1)}{\omega_c}i - c_1 = i\omega_c.
\end{equation}
Then, with definition (\ref{At}), the function $q(s)=e^{i\omega_cs} \in
C^1([-\tau_c,0],\mathbb{C})$ is an eigenvector of $\widetilde{A}_0$ associated with $i\omega_c$.

Hence, $q^* (s) = de^{i\omega_cs} \in C^1 ([0,\tau_c], \mathbb{C})$, $d\neq0$, is an eigenvector
for $\widetilde{A}_0^*$ associated with $-i\omega_c$. Moreover, we can choose $d\in\mathbb{C}$ so
that the norming condition $\langle q^*,q\rangle =1$ is satisfied. It follows that
\begin{displaymath}
\bar d = \left[1 + c_2(0)\left(\frac{\tau_c e^{-i \tau_c \omega_c}}{\omega_c}i-\frac{1-e^{-i \tau_c
\omega_c}}{\omega_c^2}\right)\right]^{-1}.
\end{displaymath}
One can check that in fact $\overline{d}=(\Delta_{\lambda}(i\omega_c,\tau_c))^{-1}$. Since
$i\omega_c$ is a simple root of $\Delta(\cdot,\tau_c)$, then $\overline{d}$ is well-defined.

From (\ref{ps2}) and (\ref{prop}) we infer also that
\begin{equation}\label{prop2}
\langle q^*,\overline{q}\rangle = 0.
\end{equation}

We are interested in the center manifold corresponding to the eigenvalue $\lambda= i\omega_c$ of
$\widetilde{A}_0$ and to the system (\ref{aeqy3}). Such a center manifold exists (see
\cite{hkw1981,c1971}): it is a locally invariant, locally attracting manifold containing the origin
and tangent at the origin to the subspace spanned by the eigenvectors corresponding to the
eigenvalues $\pm i\omega_c$ of $\widetilde{A}_0$. In fact, to reach our aim, we only need
information on the section of the center manifold, denoted $\mathcal{C}_0$, corresponding to
$\mu=0$ (see \cite{hkw1981}).

Let $y_t$ be a solution of
\begin{equation}\label{system0}
\displaystyle\frac{dy_t}{dt} = \widetilde{A}_0y_t + X_0F_0(y_t).
\end{equation}
We compute the coordinates of the section $\mathcal{C}_0$ of the center manifold corresponding to
$\mu=0$. Following the notations in \cite{hkw1981}, we define
\begin{equation}\label{z}
z(t) = \langle q^*, y_t \rangle, \qquad \textrm{ for } t\geq0.
\end{equation}
We will use $z$ and $\overline{z}$ as local coordinates of $\mathcal{C}_0$ in the directions $q^*$
and $\overline{q}^*$ respectively. We also define, for $t\geq0$ and $s\in[-\tau_c,0]$,
\begin{displaymath}
\begin{array}{rcl}
w(t,s)&=&y_t(s)-z(t)q(s)-\overline{z}(t)\overline{q}(s),\vspace{1ex}\\
&=&y_t (s) - 2\textrm{Re}[z(t)q(s)].
\end{array}
\end{displaymath}
We have
\begin{displaymath}
w(t,s) = W\left(z(t), \overline{z}(t), s\right), \quad t\geq0,\ s\in[-\tau_c,0],
\end{displaymath}
with
\begin{equation}\label{w}
W(z,\overline{z},s) = w_{20} (s) \frac{z^2}{2} + w_{11} (s) z\overline{z} + w_{02} (s)
\frac{\overline{z}^2}{2} + \ldots \end{equation}

One can notice that, for real solution $y$, $w$ is real so $w_{02} = \overline{w_{20}}$. Moreover,
(\ref{prop2}) and (\ref{z}) imply that $\langle q^*, w\rangle = 0$.

The section $\mathcal{C}_0$ of the center manifold is locally invariant under equation
(\ref{system0}): any solution that starts in it will stay in it for any time $t$ in some nontrivial
interval; therefore, if $y_t\in \mathcal{C}_0$ we have
$$
\frac{d}{dt}z(t) = \langle q^*, \widetilde{A}_0y_t + X_0F_0(y_t)\rangle
$$
so, from (\ref{z}), it follows that
\begin{equation}\label{zprime}
\setlength\arraycolsep{2pt}
\begin{array}{rcl}
\displaystyle\frac{d}{dt}z(t)&=&i\omega_c z(t) + \overline{d}
F_0\left(W(z(t),\overline{z}(t),\cdot)+2\textrm{Re}[z(t)q]\right),\vspace{1ex}\\
&=&i\omega_cz(t)+g(z(t),\overline{z}(t)),
\end{array}
\end{equation}
with
\begin{displaymath}
g(z,\overline{z})=\overline{d} F_0\left(W(z,\overline{z},\cdot)+2\textrm{Re}[zq]\right).
\end{displaymath}
We use the Taylor expansion of $f$ around $x^*$ to rewrite $F_0$ as
\begin{displaymath}
\begin{array}{rcl}
F_0(\varphi)&=&c_3\displaystyle\int^0_{-\tau_c}\left[\frac{\beta_2}{2!}\displaystyle\varphi(\theta)^2
+\frac{\beta_3}{3!}\varphi(\theta)^3
+\mathcal{O}(\varphi(\theta)^4)\right]d\theta\vspace{1ex}\\
&&-\beta_0\left[ \displaystyle\frac{\beta_2}{2!}\varphi(0)^2 +\frac{\beta_3}{3!}\varphi(0)^3+
\mathcal{O}(\varphi(0)^4)\right],
\end{array}
\end{displaymath}
where
\begin{equation}\label{parameters2}
c_3:=\displaystyle\frac{2\beta_0}{\tau_c}, \qquad \beta_2:=f^{\prime\prime}(x^*) \qquad \textrm{
and } \qquad \beta_3 := f^{\prime\prime\prime} (x^*).
\end{equation}
If we denote, for convenience, $w(s)=W(z,\overline{z},s)$, we then obtain
\begin{equation}\label{taylor}
\setlength\arraycolsep{2pt}
\begin{array}{rcl}
g(z,\overline{z}) &=& \overline{d}c_3\left\{ \displaystyle\frac{\beta_2}{2}
\displaystyle\int^0_{-\tau_c} [w(s) + ze^{i\omega_cs} + \overline{z}e^{-i\omega_cs}]^2
ds\right.\vspace{1ex}\\
&&\qquad+\displaystyle\frac{\beta_3}{6}\displaystyle\int^0_{-\tau_c}[w(s) + ze^{i\omega_cs}+\overline{z}e^{-i\omega_cs}]^3ds\vspace{1ex}\\
&&\qquad+\left.\displaystyle\int^0_{-\tau_c}\mathcal{O}([w(s)+ze^{i\omega_cs}+\overline{z}e^{-i\omega_cs}]^4)ds \right\}\vspace{1ex}\\
&-&\beta_0\overline{d}\left(\displaystyle\frac{\beta_2}{2}[w(0)+z+\overline{z}]^2 \right.\vspace{1ex}\\
&&\qquad+\left.\displaystyle\frac{\beta_3}{6}[w(0)+z+\overline{z}]^3+\mathcal{O}([w(0)+z+\overline{z}]^4)\right).
\end{array}
\end{equation}

Equation (\ref{zprime}) is called the normal form obtained by the restriction of the flow to the
center manifold. Our next goal is to compute some coefficients in the Taylor series of $g$ and to
use them to study stability of the limit cycle by computing also the Lyapunov coefficient. This
latter is given by some coefficients in the Taylor expansion  of $g(z,\overline{z})$ given by
(\ref{taylor}).
 This means, in fact, that stability of the limit cycle  is
determined by the normal form obtained through the restriction of the flow to the center manifold.

Restricting the Taylor expansion in (\ref{taylor}) to terms of order less or equal to three, we get
\begin{equation}\label{g}
g(z,\overline{z})=\frac{1}{2}g_{20}z^2+g_{11}z\overline{z}+\frac{1}{2}g_{02}{\overline{z}}^2
+\frac{1}{2}g_{21}z^2\overline{z}+\ldots \end{equation} with
\begin{equation}\label{gparameters}
\setlength\arraycolsep{1pt}
\begin{array}{rcl}
g_{20}&=&-\bar d \beta_2
\left(\beta_0+\displaystyle\frac{c_3(1-e^{-2i\omega_c\tau_c})}{2\omega_c}i\right),\vspace{1ex}\\
g_{11}&=&\bar d \beta_2(c_3 \tau_c - \beta_0),\vspace{1ex}\\
g_{02}&=&-\bar d \beta_2\left(\beta_0-\displaystyle\frac{c_3(1-e^{2i\omega_c\tau_c})}{2\omega_c}i\right),\vspace{1ex}\\
g_{21}&=&\overline{d}\bigg( c_3\left\{ \beta_2
\!\displaystyle\int^0_{-\tau_c}\!\![w_{20}(s)e^{-i\omega_cs}+2w_{11}(s)e^{i\omega_cs}]ds-\beta_3
\displaystyle\frac{(1-e^{-i\omega_c\tau_c})}{\omega_c}i\right\}\vspace{1ex}\\
&&\qquad-\beta_0\beta_2[w_{20}(0)+2w_{11}(0)]-\beta_0\beta_3\bigg).
\end{array}
\end{equation}
With the definition of $c_3$ given in (\ref{parameters2}), one can see that
\begin{displaymath}
g_{11}=\overline{d}\beta_0\beta_2.
\end{displaymath}
One can note that we do not need the coefficients $g_{ij}$ with $i+j>2$ except $g_{21}$ (see
\cite{hkw1981}).

At this point, we still need to compute $w_{20}(s)$ and $w_{11}(s)$ for $s\in[-\tau_c,0]$. It
follows directly from (\ref{system0}), (\ref{z}) and (\ref{zprime}), that
\begin{equation}\label{wprime}
\begin{array}{rcl}
\displaystyle\frac{d}{dt}w(t,\cdot)&=&\displaystyle\frac{d}{dt}y_t-
\frac{d}{dt}\left[z(t)q-\overline{z}(t)\overline{q}\right],\vspace{1ex}\\
&=&\widetilde{A}_0w(t,\cdot)+X_0F_0\left(w(t,\cdot)+2\textrm{Re}(z(t)q)\right)\vspace{1ex}\\
&&-2\textrm{Re}\left(g(z(t),\overline{z}(t))q\right),\vspace{1ex}\\
&=&\widetilde{A}_0w(t,\cdot)+H(z(t),\overline{z}(t),\cdot),
\end{array}
\end{equation}
where, for $s\in[-\tau_c,0]$,
\begin{equation}\label{H}
H(z,\overline{z},s)= -2\textrm{Re}\left(g(z,\overline{z})q(s)\right)
+X_0(s)F_0\left(W(z,\overline{z},\cdot)+2\textrm{Re}(zq)\right).
\end{equation}
Let $s\in[-\tau_c,0)$ be fixed. From (\ref{H}) with (\ref{g}) we have
\begin{eqnarray}
H(z,\overline{z},s)&=&-2\textrm{Re}\left(g(z,\overline{z})q(s)\right),\nonumber\\
&=&-g(z,\overline{z})q(s)-\bar g(z,\overline{z})\bar q(s),\nonumber\\
&=&-\left(g_{20}\frac{z^2}{2}+g_{11}z\overline{z}+g_{02}\frac{{\overline{z}}^2}{2}+\ldots\right)q(s)\nonumber\\
&&-\left(\overline{g}_{20}\frac{{\overline{z}}^2}{2}+\overline{g}_{11}z\overline{z}+\overline{g}_{02}\frac{z^2}{2}+\ldots\right)\overline{q}(s).
\label{H2}
\end{eqnarray}
Considering the expansion
$$
H(z,\overline{z},s) = H_{20}(s)\frac{z^2}{2} + H_{11}(s)
z\overline{z}+H_{02}(s)\frac{\overline{z}^2}{2} + \ldots $$ and comparing the coefficients with
those in (\ref{H2}) yields
\begin{equation}\label{H3}
\begin{array}{lll}
H_{20}(s)=-g_{20}q(s)-\overline{g}_{02}\overline{q}(s),\vspace{1ex}\\
H_{11}(s)=-g_{11}q(s)-\overline{g}_{11}\overline{q}(s),\vspace{1ex}\\
H_{02}(s)=\overline{H}_{20}(s).
\end{array}
\end{equation}
From (\ref{w}), (\ref{g}) and (\ref{wprime}), we obtain
\begin{equation}\label{H4}
\begin{array}{l}
(\widetilde{A}_0-2i\omega_c)w_{20}(s)=-H_{20}(s),\vspace{1ex}\\
\widetilde{A}_0w_{11}(s)=-H_{11}(s),\vspace{1ex}\\
(\widetilde{A}_0+2i\omega_c)w_{02}(s)=-H_{02}(s).
\end{array}
\end{equation}
Using (\ref{At}) and (\ref{H3}), (\ref{H4}) becomes
\begin{displaymath}
\begin{array}{ll}
\displaystyle\frac{dw_{20}}{ds}(s) = 2i\omega_cw_{20}(s)+g_{20}q(s)+\overline{g}_{02}\overline{q}(s)\vspace{1ex},\\
\displaystyle\frac{dw_{11}}{ds}(s) = g_{11}q(s)+\overline{g}_{11}\overline{q}(s).
\end{array}
\end{displaymath}
Solving the above system, we obtain:
\begin{equation}\label{w20}
w_{20}(s)=-\frac{g_{20}}{i\omega_c}e^{i\omega_cs}-\frac{\overline{g}_{02}}{3i\omega_c}e^{-i\omega_cs}+E_1e^{2i\omega_cs},
\end{equation}
\begin{equation}\label{w11}
w_{11}(s)=\frac{g_{11}}{i\omega_c}e^{i\omega_cs}-\frac{\overline{g}_{11}}{i\omega_c}e^{-i\omega_cs}+E_2,
\end{equation}
where $E_1$ and $E_2$ can be determined by setting $s=0$ in (\ref{H}). In fact, we have
$$
H(z,\overline{z},0)=-2\textrm{Re}\left(g(z,\overline{z})\right)
+F_0\left(W(z,\overline{z},\cdot)+2\textrm{Re}(zq)\right),
$$
so we deduce
\begin{displaymath}
\begin{array}{rcl}
H_{20}(0)&=&-g_{20}-\overline{g}_{02}+\beta_2\left(c_3\displaystyle\int^0_{-\tau_c}e^{2i\omega_cs}ds-\beta_0\right),\vspace{1ex}\\
H_{11}(0)&=&-g_{11}-\overline{g}_{11}+\beta_2\left(c_3\displaystyle\int^0_{-\tau_c}ds-\beta_0\right).
\end{array}
\end{displaymath}
From (\ref{At}) and (\ref{H4}), we have:
\begin{eqnarray}
c_2(0)\int^0_{-\tau_c}w_{20}(s)ds-(c_1+2i\omega_c)w_{20}(0)&=&-H_{20}(0),\label{H6}\\
c_2(0)\int^0_{-\tau_c}w_{11}(s)ds-c_1w_{11}(0)&=&-H_{11}(0).\label{H7}
\end{eqnarray}
Substituting (\ref{w20}) in (\ref{H6}) and using (\ref{prop}) we eventually get
\begin{equation}\label{E1}
E_1=-\frac{\beta_2(c_3(1-e^{-2i\omega_c\tau_c})-2i\omega_c\beta_0)}{c_2(0)(1-e^{-2i\omega_c\tau_c})-2c_1i\omega_c+4\omega_c^2}.
\end{equation}
Similarly, substituting (\ref{w11}) in (\ref{H7}) we get
\begin{equation}\label{E2}
E_2=-\frac{\beta_2(c_3\tau_c-\beta_0)}{c_2(0)\tau_c-c_1}.
\end{equation}
Using the definitions given in (\ref{parameters1}) and (\ref{parameters2}), we can write
\begin{displaymath}
E_2=\frac{\beta_0\beta_2}{\delta-\beta_0\beta_1}.
\end{displaymath}
We are now able to complete the calculation of $g_{21}$ in (\ref{gparameters}) using the above
values of $w_{20}$, $w_{11}$, $E_1$ and $E_2$. We set
\begin{displaymath}
K := \int^0_{-\tau_c}\left[w_{20}(s)e^{-i\omega_cs}+2w_{11}(s)e^{i\omega_cs}\right]ds.
\end{displaymath}
Then, from (\ref{w20}) and (\ref{w11}), we obtain
\begin{equation}\label{K}
K=\displaystyle\frac{(\overline{g}_{02}+6g_{11})(1-e^{-2i\omega_c\tau_c})}{6\omega_c^2}\vspace{1ex}\\
+\frac{\tau_c(g_{20}+2\overline{g}_{11})-(E_1+2E_2)(1-e^{-i\omega_c\tau_c})}{\omega_c}i
\end{equation}
where $E_1$ and $E_2$ are given respectively by (\ref{E1}) and (\ref{E2}).

From (\ref{gparameters}), (\ref{w20}) and (\ref{w11}), we have
\begin{displaymath}
\begin{array}{rcl}
g_{21}&=&\overline{d}\bigg( c_3\beta_2K-c_3\beta_3
\displaystyle\frac{(1-e^{-i\omega_c\tau_c})}{\omega_c}i\vspace{1ex}\\
&&-\beta_0\beta_2\left(E_1+2E_2+\displaystyle\frac{3g_{20}+\overline{g}_{02}}{3\omega_c}i
+\displaystyle\frac{2(\overline{g}_{11}-g_{11})}{\omega_c}i\right)-\beta_0\beta_3\bigg),
\end{array}
\end{displaymath}
where $K$ is given by (\ref{K}).

Based on the above analysis and calculation, we can see that each $g_{ij}$ in (\ref{gparameters})
is determined by the parameters and delay in equation (\ref{eq}). Thus we can explicitly compute
the following quantities:
\begin{equation}\label{quantities}
\begin{array}{rcl}
L_1(0)&=&\displaystyle\frac{i}{2\omega_c}\left(g_{20}g_{11}-2|g_{11}|^2-\frac{1}{3}|g_{02}|^2\right)+\frac{1}{2}g_{21},\vspace{1ex}\\
l_1(0)&=&\textrm{Re}\left(L_1(0)\right),\vspace{1ex}\\
\mu_2&=&-\displaystyle\frac{l_1(0)}{\textrm{Re}(\lambda^{'}(\tau_c))},\vspace{1ex}\\
b_2&=&2l_1(0),\vspace{1ex}\\
T_2&=&-\displaystyle\frac{\textrm{Im}(L_1(0))+\mu_2\textrm{Im}\lambda^{'}(\tau_c)}{\omega_c}.
\end{array}
\end{equation}
One knows from \cite{hkw1981} that the following properties hold: if $\mu_2>0$ ($<0$) then the Hopf
bifurcation is supercritical (subcritical) and the bifurcating periodic solutions exist for
$\tau>\tau_c$ ($\tau<\tau_c$); solutions are orbitally stable (unstable) if $b_2<0$ ($>0$); and the
period of bifurcating periodic solution increases (decreases) if $T_2>0$ ($<0$).

The coefficient $\lambda^{\prime}(\tau_c)$ in (\ref{quantities}) is given by (\ref{muprimezero})
and (\ref{omegaprimezero}). In particularly, we have proved in Section \ref{shopf}, property
(\ref{muprime}), that
\begin{displaymath}
\textrm{Re}(\lambda^{\prime}(\tau_c))>0.
\end{displaymath}
In summary, this leads to the following result:

\begin{theorem}
If the Lyapunov coefficient $l_1(0)$, defined in (\ref{quantities}), is negative (resp. positive)
then the Hopf bifurcation is supercritical (subcritical) and the bifurcating periodic solutions
exist for $\tau>\tau_c$ ($\tau<\tau_c$), and solutions are orbitally stable (unstable); The
coefficient $T_2$ determines the period of the bifurcating periodic solutions: the period increases
(decreases) if $T_2>0$ ($<0$).
\end{theorem}

\section{Numerical results and simulations}\label{snumerical}

We numerically compute, in this section, the formulas obtained above to determine the behavior of
the periodic solutions of equation (\ref{eq}).

We choose $f$ as in (\ref{functionf}). In order to satisfy the assumptions of Theorem
\ref{theoremhopfbifurcation}, we have to choose $\delta$, $\beta_0$ and $n$ such that $\beta_1<0$
and (\ref{hyph}) holds true; that is
\begin{displaymath}
n>\frac{2(1-h(x_0))}{1-2h(x_0)}\ \frac{\beta_0}{\beta_0-\delta}.
\end{displaymath}
We take $\delta$ and $\beta_0$ as given in (\ref{parametersvalues}). Then the above conditions are
in fact satisfied for $n\geq 2.42$.

Using Maple 9, we are able to compute the coefficients in (\ref{quantities}), listed in the
following table for $n\in\{3,4,5,6,7,8\}$:

\smallskip

\begin{center}
{\small \begin{tabular}{|c|c|c|c|c|c|}
\hline &&&&&\\
{\bf $n$} & {\bf $x^*$} & {\bf $\tau_c$} & {\bf $\omega_c$} & {\bf $l_1(0)$} & {\bf $T_2$}\\
&&&&&\\
\cline{1-6} $3$ & 3.2523 & 18.1270 & 0.1380 & -0.0026 & 0.1540  \\
\cline{1-6} $4$ & 2.4218 & 10.9688 & 0.2091 & -0.0375 & 0.8083  \\
\cline{1-6} $5$ & 2.0291 & 7.8748  & 0.2785 & -0.1504 & 2.2097  \\
\cline{1-6} $6$ & 1.8034 & 6.1447  & 0.3472 & -0.3975 & 4.5441  \\
\cline{1-6} $7$ & 1.6577 & 5.0385  & 0.4155 & -0.8399 & 7.9335  \\
\cline{1-6} $8$ & 1.5562 & 4.2702  & 0.4836 & -1.5412 & 12.4589 \\
\hline
\end{tabular}}
\end{center}


Even though we do not give values of $l_1(0)$ and $T_2$ for all $n\geq 2.42$, we can notice that
observations indicate that $l_1(0)$ is strictly negative and $T_2$ strictly positive for $n\geq
2.42$. Hence the unique Hopf bifurcation of equation (\ref{eq}) seems to be supercritical and
solutions orbitally stable, with increasing periods.

\begin{figure}[!hpt]
\begin{center}R
\includegraphics[width=6.5cm, height=4cm]{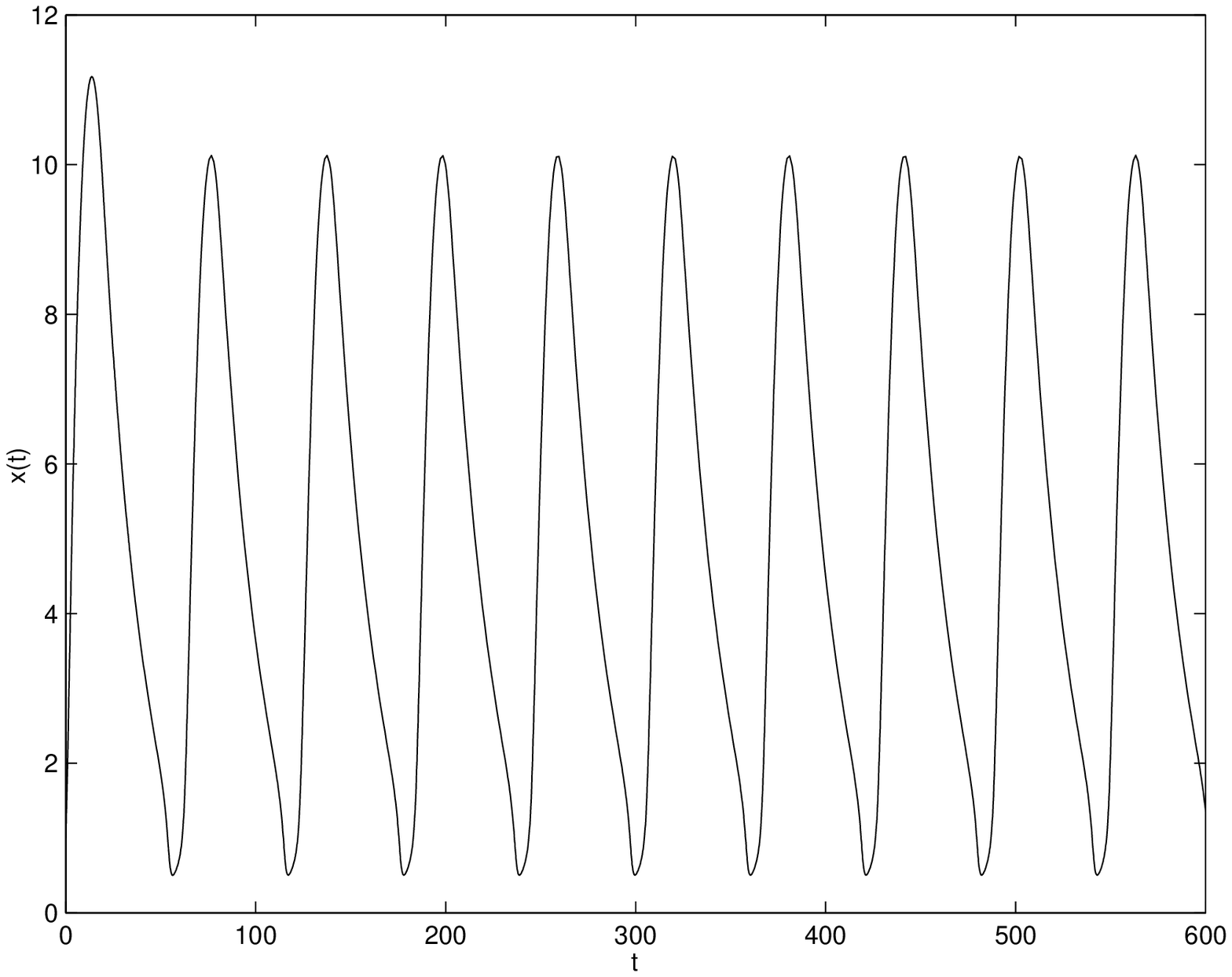}
\includegraphics[width=6.5cm, height=4cm]{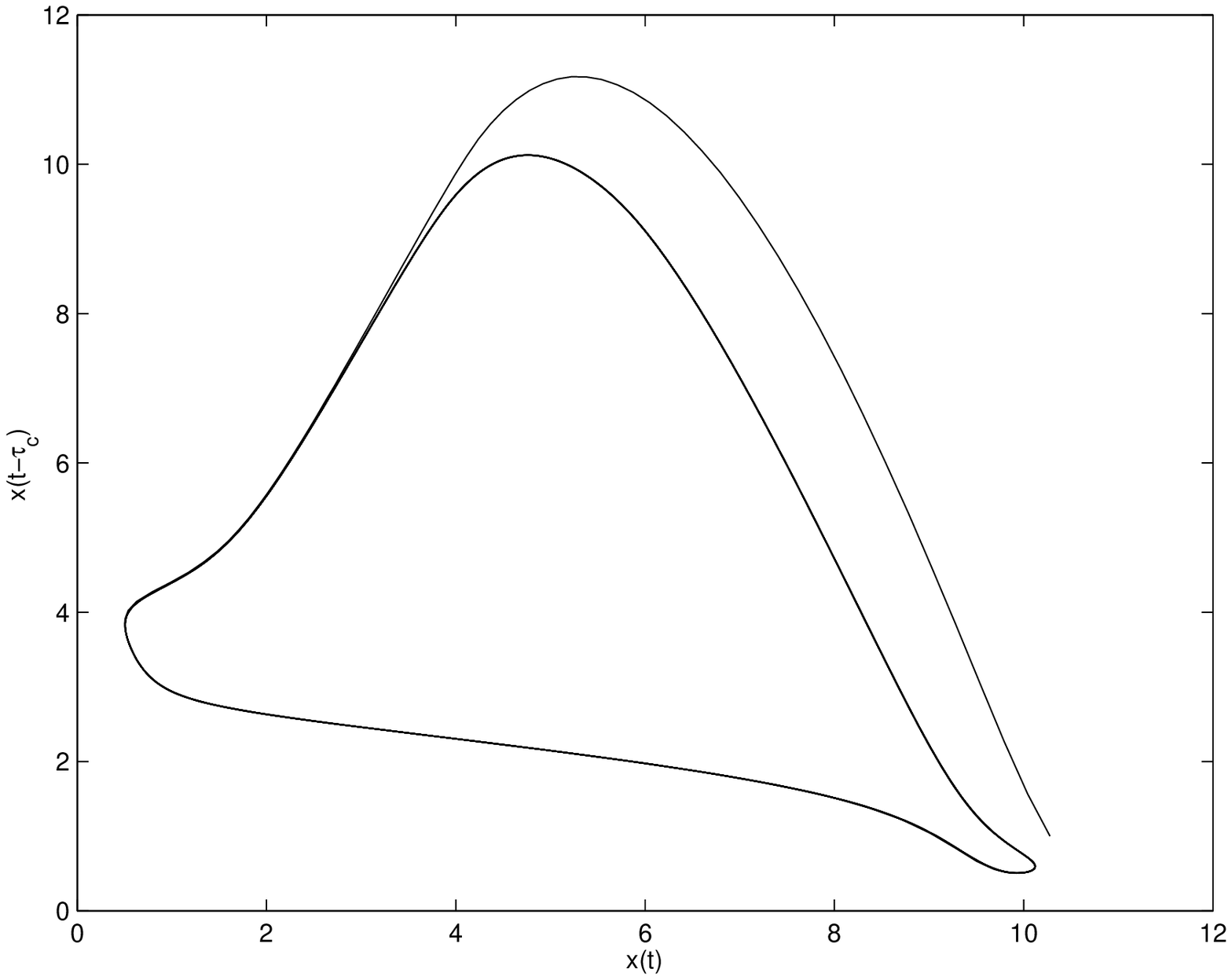}
\end{center}
\caption{Solutions of (\ref{eq}) are displayed when $n=3$. Periods of the oscillations are close to
46 days.}\label{fign3}
\end{figure}
\begin{figure}[!hpt]
\begin{center}
\includegraphics[width=6.5cm, height=4cm]{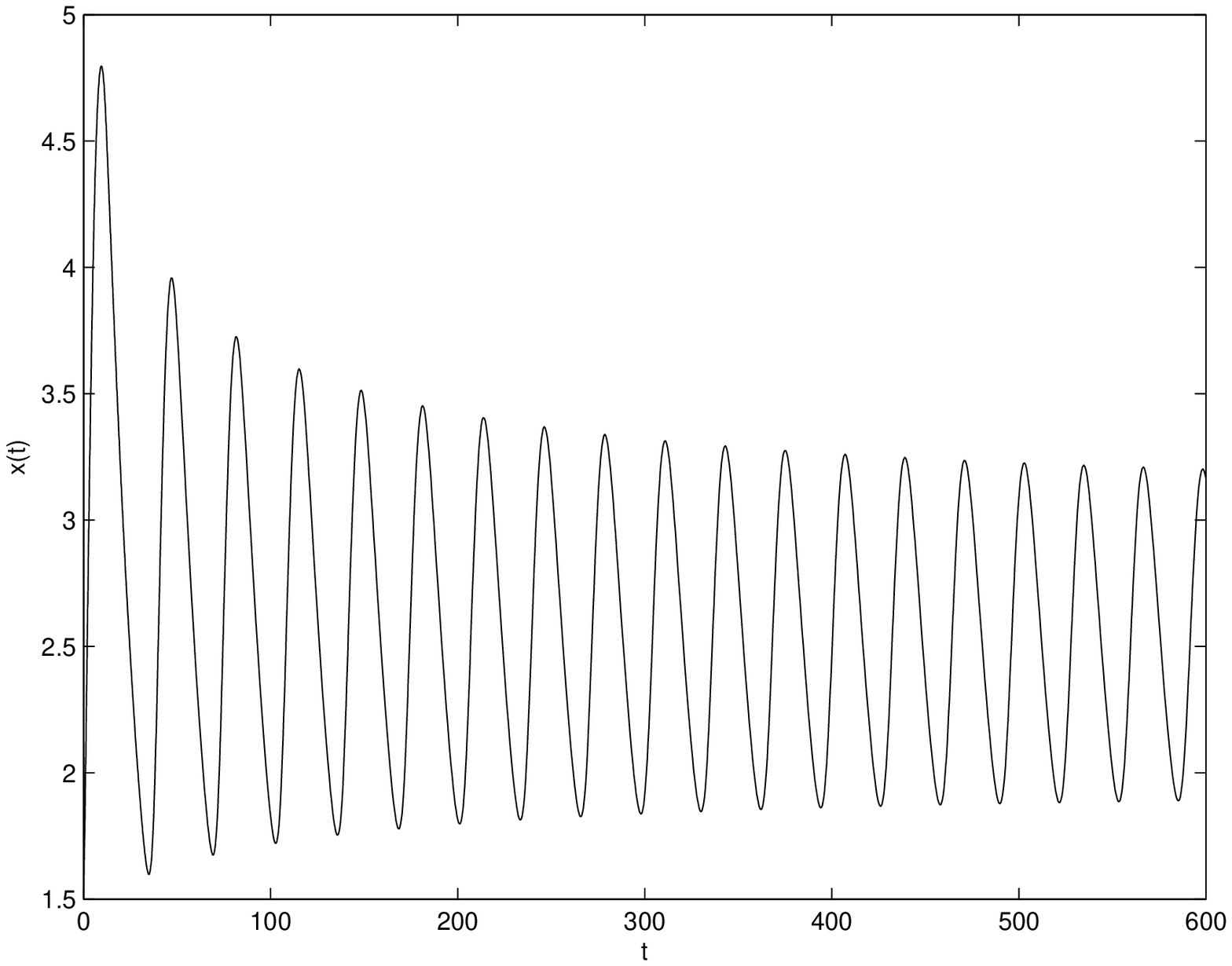}
\includegraphics[width=6.5cm, height=4cm]{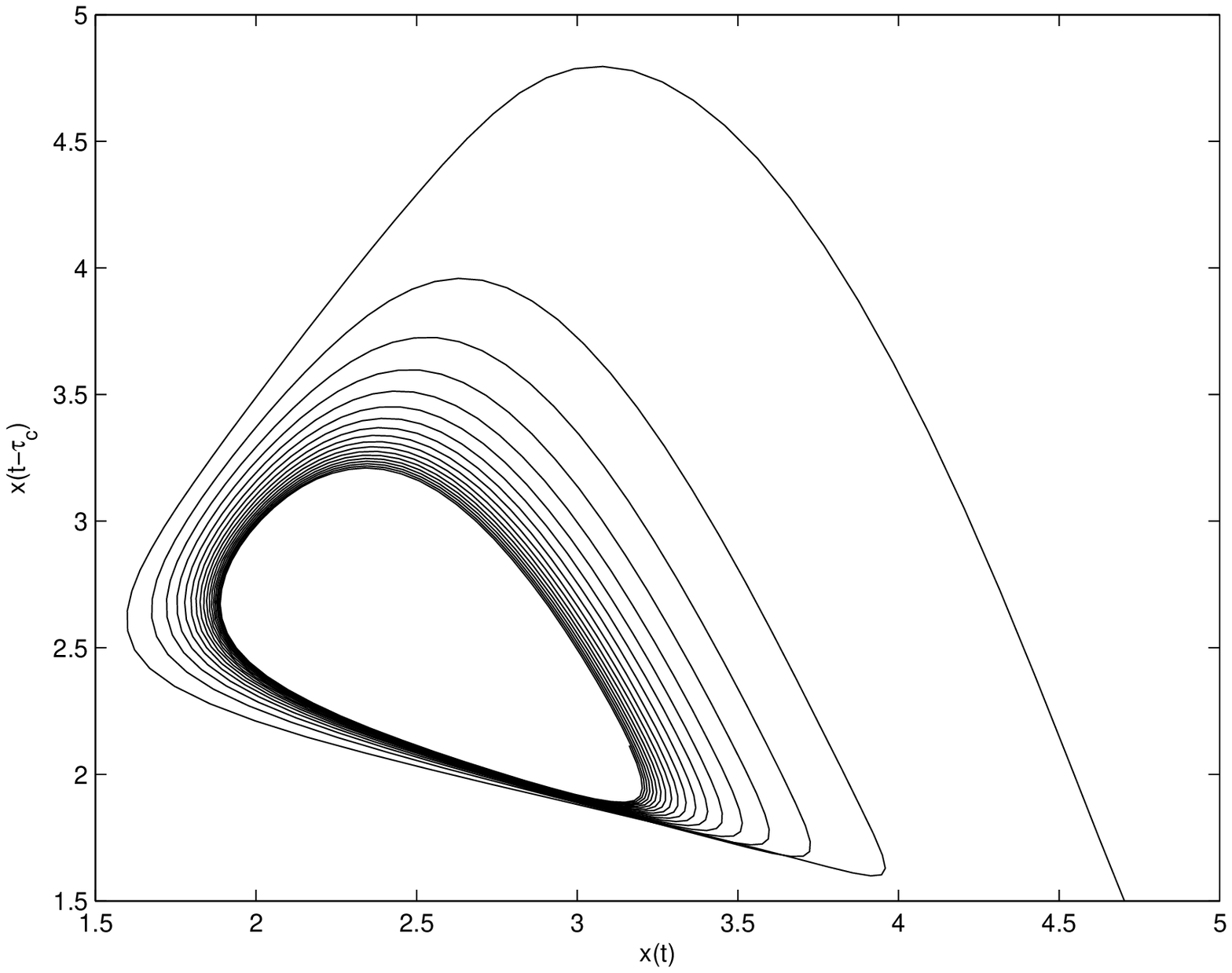}
\end{center}
\caption{Solutions of (\ref{eq}) are displayed when $n=4$. Periods of the oscillations are close to
30 days.}\label{fign4}
\end{figure}
\begin{figure}[!hpt]
\begin{center}
\includegraphics[width=6.5cm, height=4cm]{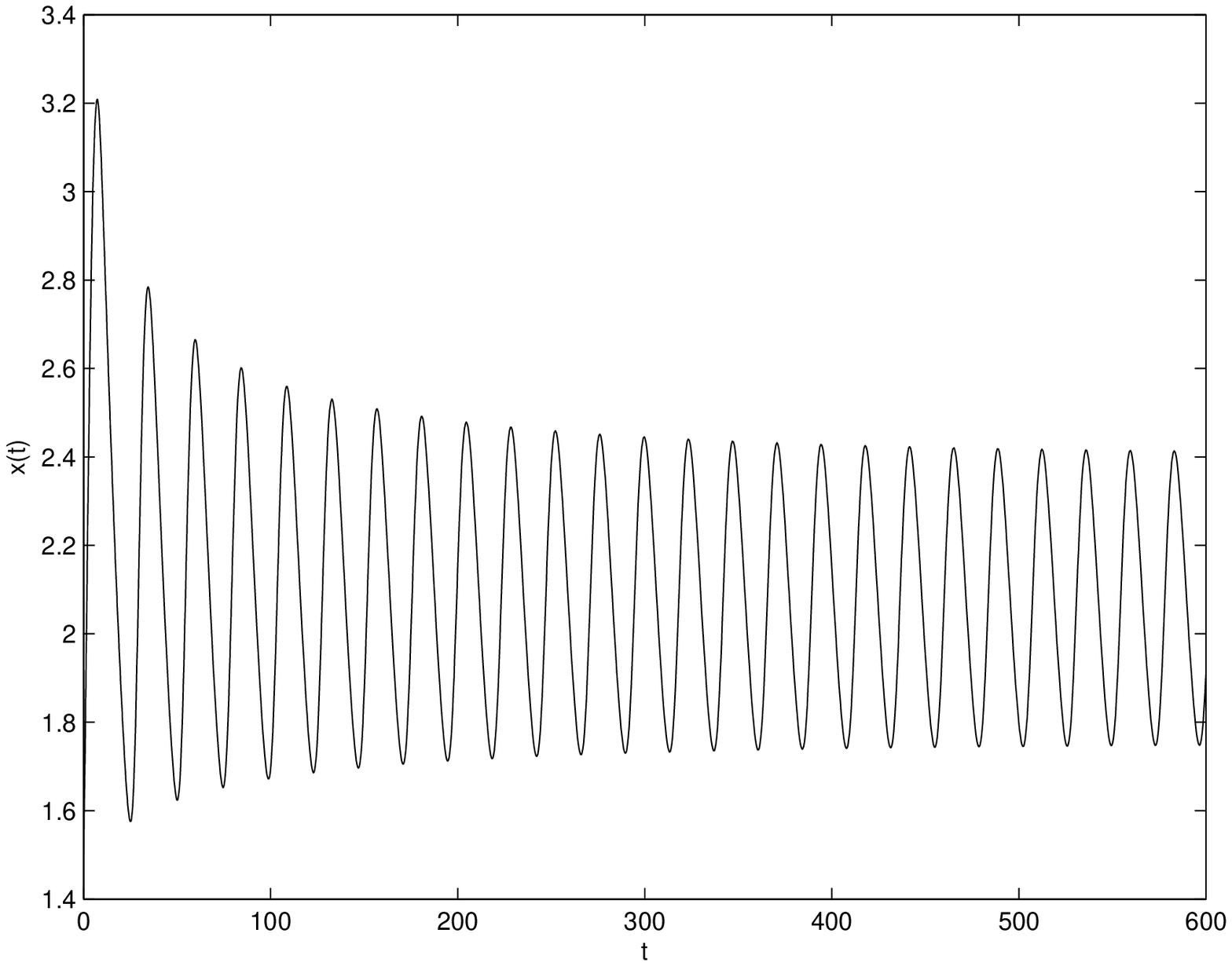}
\includegraphics[width=6.5cm, height=4cm]{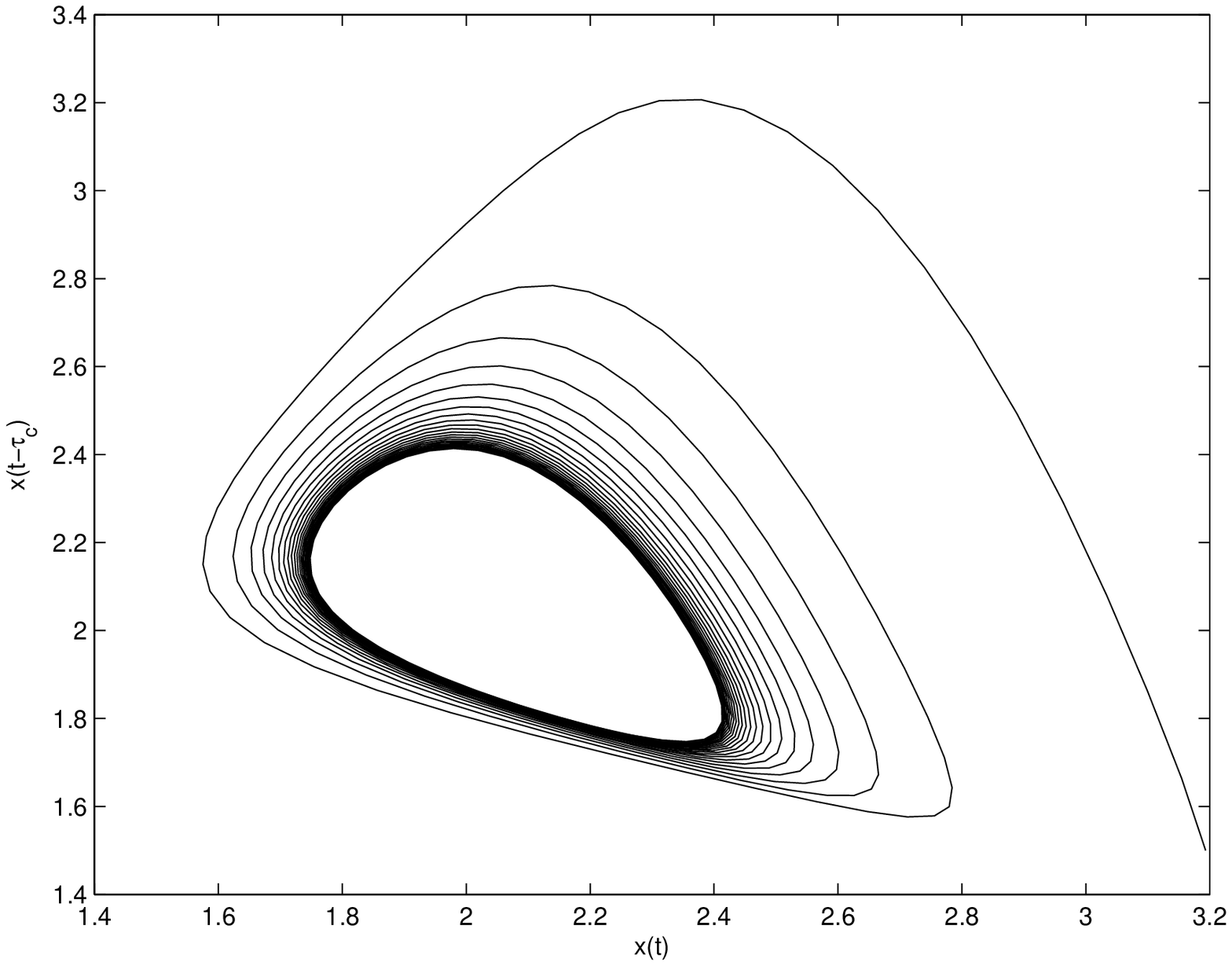}
\end{center}
\caption{Solutions of (\ref{eq}) are displayed when $n=5$. Periods of the oscillations are close to
23 days.}\label{fign5}
\end{figure}
\begin{figure}[!hpt]
\begin{center}
\includegraphics[width=6.5cm, height=4cm]{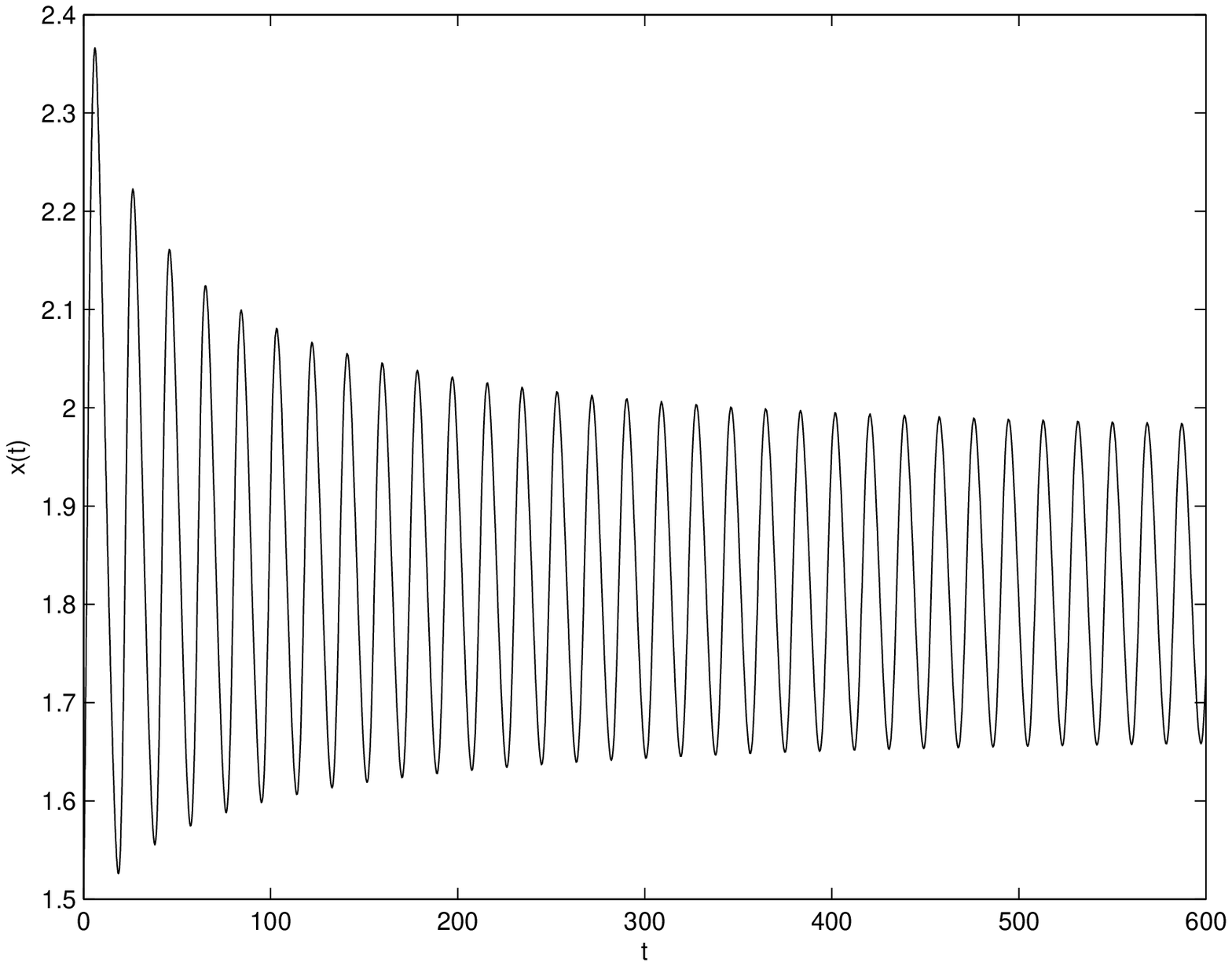}
\includegraphics[width=6.5cm, height=4cm]{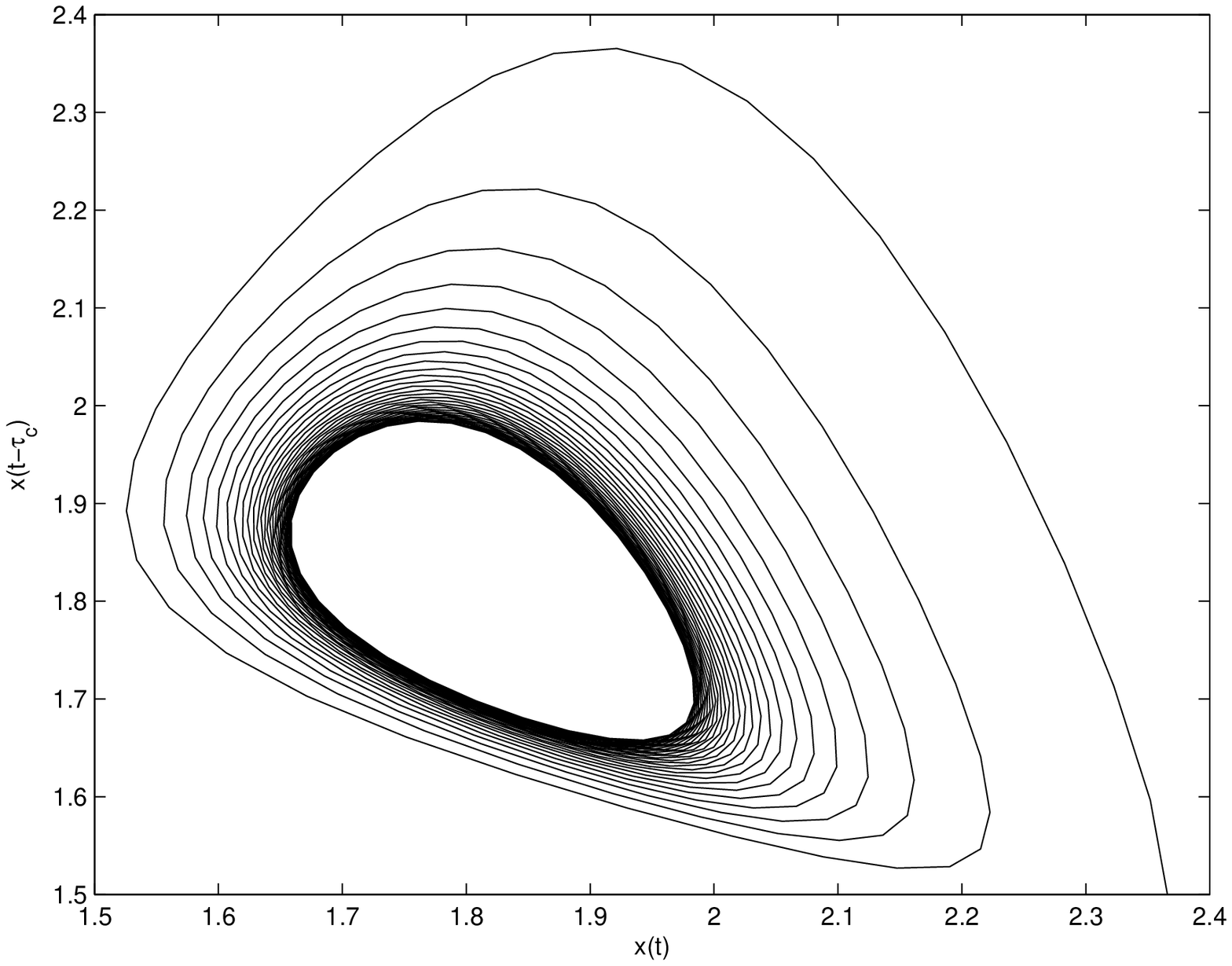}
\end{center}
\caption{Solutions of (\ref{eq}) are displayed when $n=6$. Periods of the oscillations are close to
18 days.}\label{fign6}
\end{figure}
\begin{figure}[!hpt]
\begin{center}
\includegraphics[width=6.5cm, height=4cm]{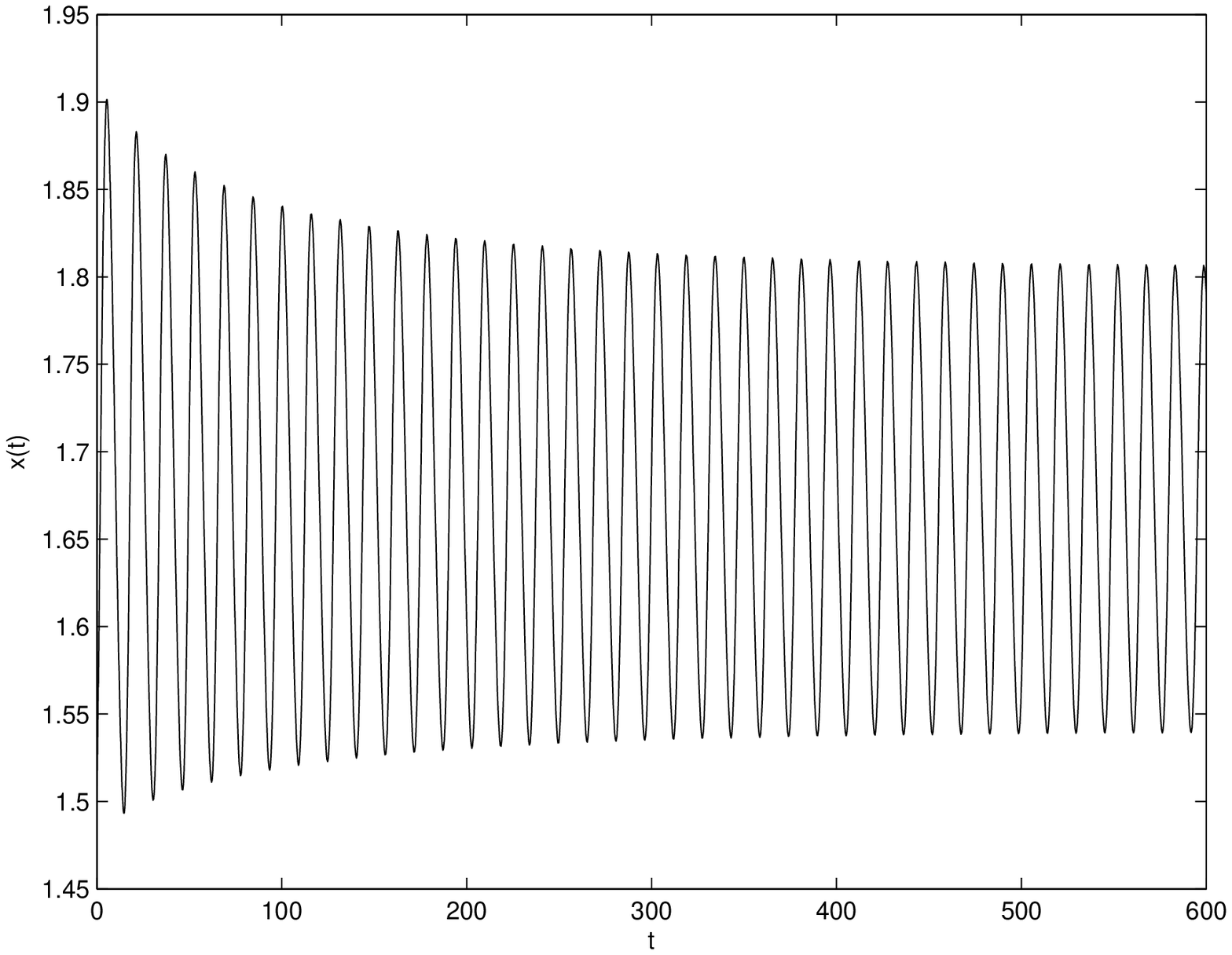}
\includegraphics[width=6.5cm, height=4cm]{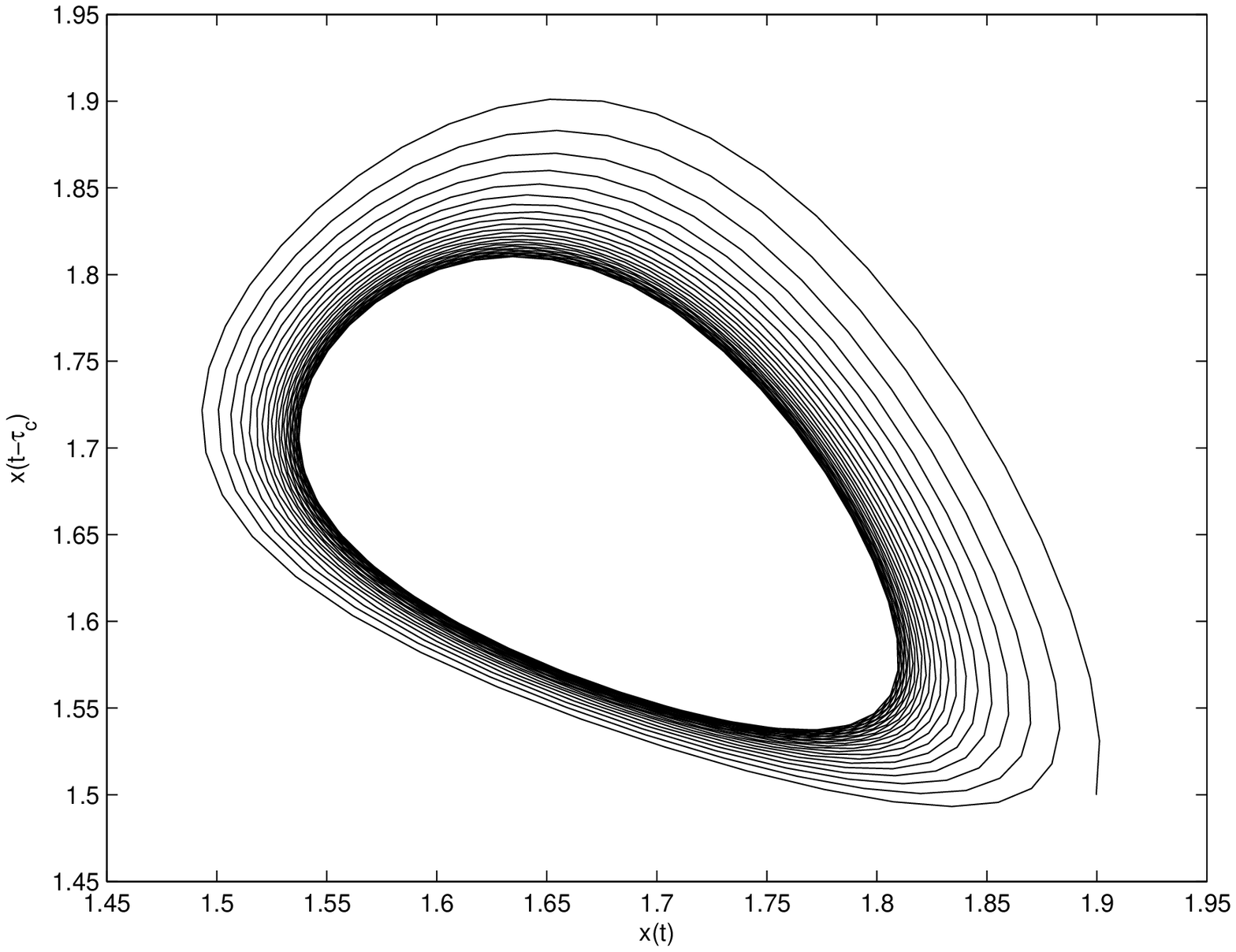}
\end{center}
\caption{Solutions of (\ref{eq}) are displayed when $n=7$. Periods of the oscillations are close to
15 days.}\label{fign7}
\end{figure}
\begin{figure}[!hpt]
\begin{center}
\includegraphics[width=6.5cm, height=4cm]{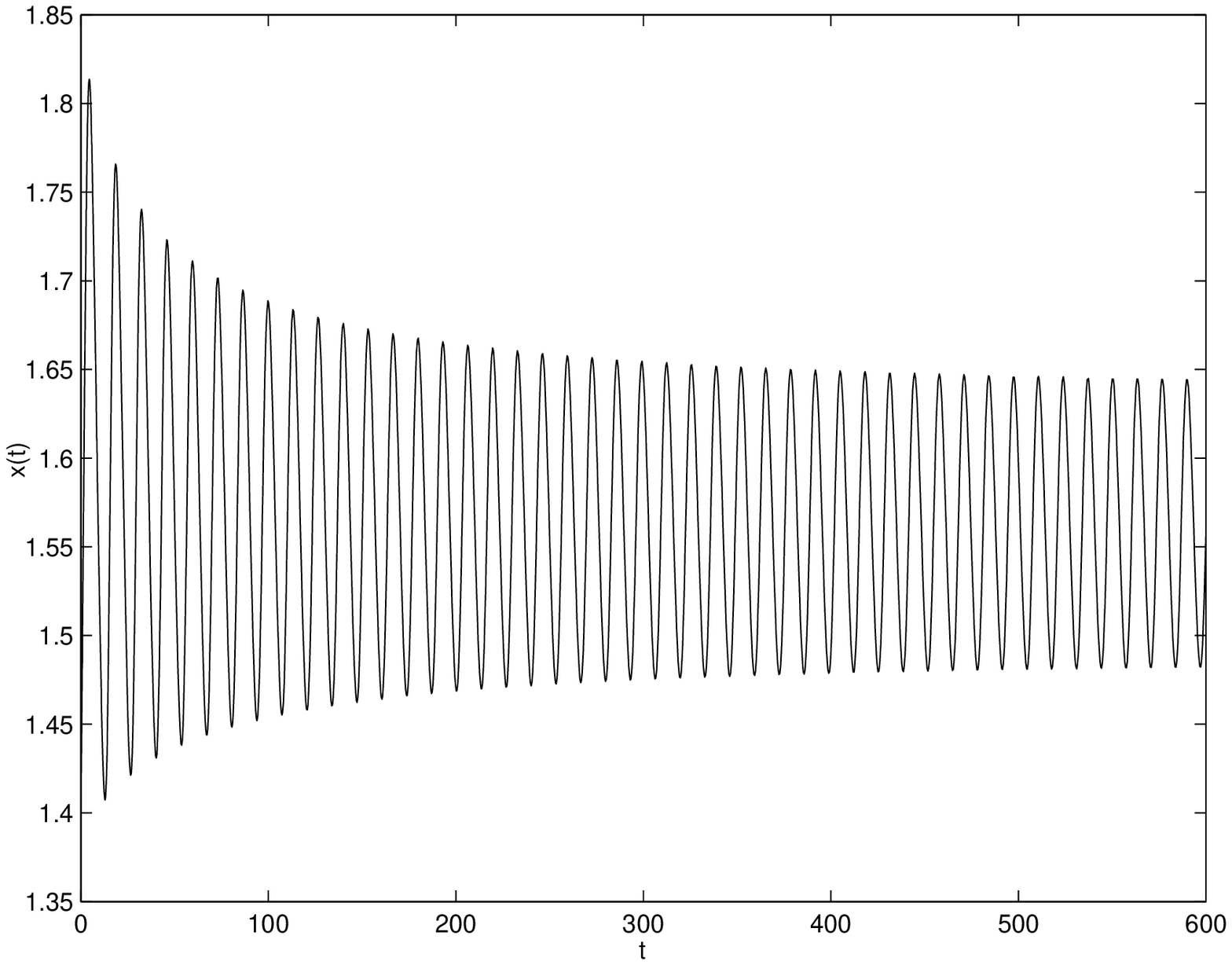}
\includegraphics[width=6.5cm, height=4cm]{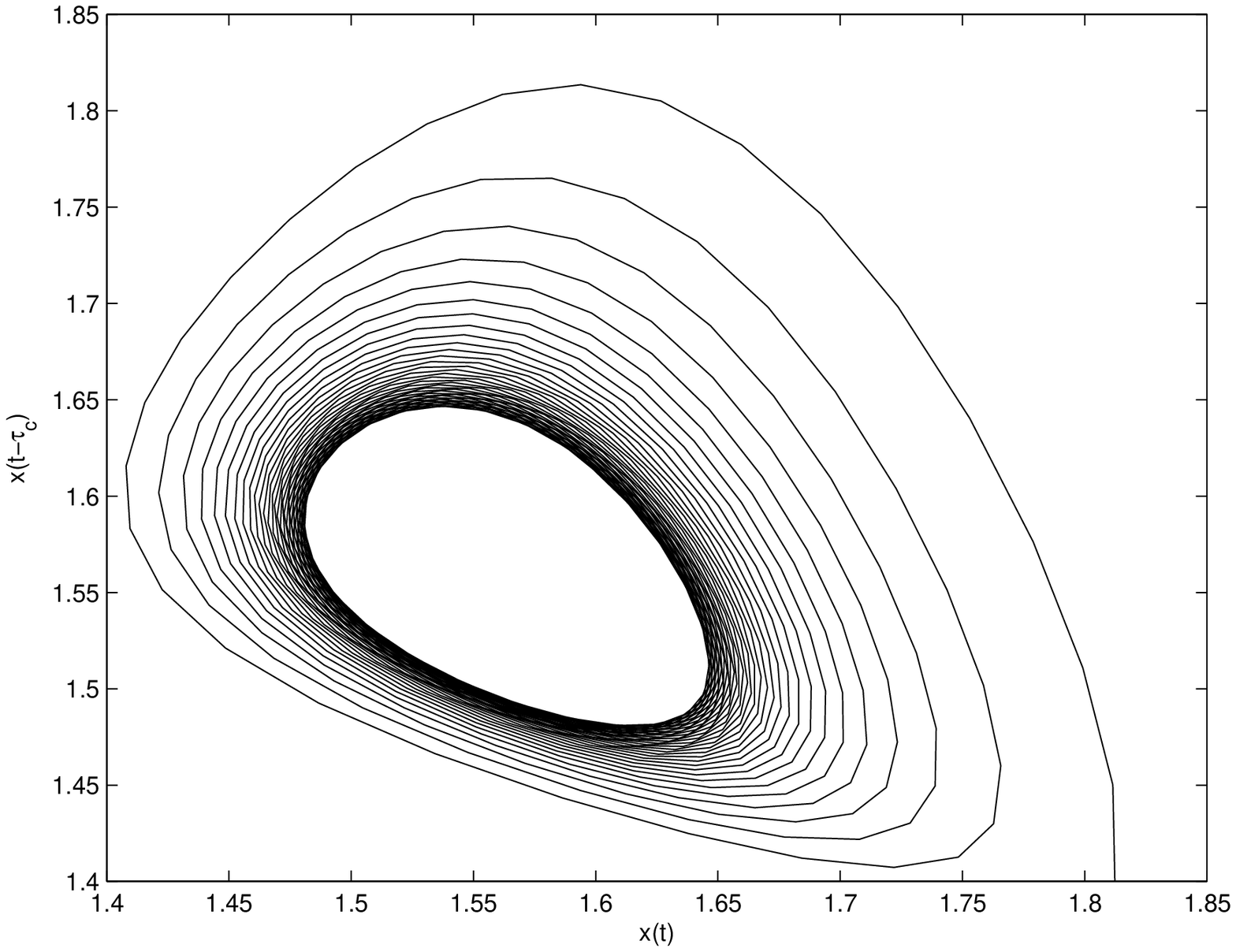}
\end{center}
\caption{Solutions of (\ref{eq}) are displayed when $n=8$. Periods of the oscillations are close to
13 days.}\label{fign8}
\end{figure}

Using the {\sc Matlab} solver dde23 \cite{dde23}, we can compute the solutions of equation
(\ref{eq}) for the above-mentioned values of the parameter $n$ and for any positive initial
condition. Solutions of (\ref{eq}) are shown versus time $t$ and in the phase plane in Fig
\ref{fign3} to \ref{fign8}.


\section{Discussion}\label{sdiscussion}

Many hematological diseases involve oscillations about a steady-state during the chronic period.
These oscillations give rise to instability in the hematopoietic stem cell count. Chronic
myelogenous leukemia (see Fortin and Mackey \cite{fm1999}) is one of the most common types of
hematological disease characterized by the existence of periodic oscillations (oscillations of
leukocytes with periods from 30 to 100 days). Experimental observations have led to the conclusion
that this dynamic instability is located in the hematopoietic stem cells compartment.

We have studied, in this paper, a mathematical model of pluripotent hematopoietic stem cells
dynamics in which the length of the proliferating phase is uniformly distributed on an interval. We
have shown that instability can occur in this model via a Hopf bifurcation, leading to periodic
solutions usually orbitally stable with increasing periods. This has been obtained throughout the
description of a center manifold and the subsequently study of the normal form.

Periods of the oscillations obtained in numerical simulations, in Section \ref{snumerical}, may be
in the order of 30 to 50 days (at the bifurcation) when the parameter $n$ is not too large,
corresponding to what can be observed with chronic myelogenous leukemia. It has already been
noticed by Pujo-Menjouet and Mackey \cite{pm2004} that this parameter $n$, which describes the
sensitivity of the rate of reintroduction $\beta$, plays a crucial role in the appearance of
periodic solutions when the delay is constant. The sensitivity $n$ describes the way the rate of
introduction in the proliferating phase reacts to changes in the resting phase population produced
by external stimuli: a release of erythropoietin, for example, or the action of some growth
factors. Since periodic hematological diseases are supposed to be due to hormonal control
destabilization (see \cite{fm1999}), then $n$ seems to be appropriate to identify causes leading to
periodic solutions.

\end{document}